\documentclass{article}

\usepackage{verbatim} \usepackage{geometry} \usepackage{amsmath}
\usepackage{amssymb} \usepackage{amsfonts} \usepackage{mathtools}
\usepackage{amsthm}
\usepackage{mathabx}
\usepackage{listings} \usepackage{graphicx}
\usepackage{color} \usepackage{stmaryrd} \usepackage{cellspace}
\usepackage{makecell}
\usepackage{inputenc}
\usepackage[colorinlistoftodos,obeyDraft]{todonotes}
\usepackage{dsfont}
\usepackage{ebproof}

\usepackage{tikz}

\usetikzlibrary{positioning,calc,backgrounds}

\usepackage{quiver}
\usepackage[ backend=biber, style=numeric, citestyle=numeric,
maxbibnames=99, isbn=false, doi=false, eprint=false, block=ragged,
uniquename=init ]{biblatex}

\DefineBibliographyStrings{english}{%
  bibliography = {References}, }

\usepackage{hyperref} \usepackage[nameinlink,capitalise]{cleveref}
\hypersetup{final}

\newtheorem{theorem}{Theorem} \newtheorem{prop}[theorem]{Proposition}
\newtheorem{cor}[theorem]{Corollary}
\newtheorem{lemma}[theorem]{Lemma} \theoremstyle{definition}
\newtheorem{definition}[theorem]{Definition} \theoremstyle{remark}

\newtheorem{example}{Example}

\DeclareMathOperator{\id}{id}
\newcommand*{\Coh}[3]{\ensuremath\mathsf{coh}\;(#1:#2)[#3]}
\newcommand*{\Ctx}{\ensuremath{\mathsf{Ctx}}}
\newcommand*{\Sub}{\ensuremath{\mathsf{Sub}}}
\newcommand*{\Type}{\ensuremath{\mathsf{Type}}}
\newcommand*{\Term}{\ensuremath{\mathsf{Term}}}
\newcommand*{\Var}{\ensuremath{\mathsf{Var}}}
\newcommand*{\Catt}{\ensuremath{\mathsf{Catt}}}
\newcommand*{\Cattsa}{\ensuremath{\mathsf{Catt_{sa}}}}
\newcommand*{\arr}[3]{\ensuremath{#1 \to_{#2} #3}}
\newcommand*{\FV}{\ensuremath{\mathsf{FV}}}
\newcommand*{\sub}[2]{\ensuremath{#1\llbracket #2 \rrbracket}}
\newcommand*{\drop}{\ensuremath{\mathsf{drop}}}
\newcommand*{\supp}{\ensuremath{\mathsf{supp}}}

\renewcommand*{\L}{\ensuremath{\mathsf{L}}}
\newcommand*{\insertion}[3]{\ensuremath{#1\ll\!\!#2\ #3}}

\newcommand{\LM}{\ensuremath{\mathsf{LM}}}
\newcommand*{\rto}{\ensuremath{\rightsquigarrow}}
\newcommand\doubleplus{+\kern-1.3ex+\kern0.8ex}
\newcommand*{\sd}{\ensuremath{\mathsf{sd}}}

\newcommand*{\+}{\mathbin{\hash}}
\DeclareMathOperator*{\bighash}{\text{\LARGE \(\+\)}}

\newcolumntype{M}[1]{>{\centering\arraybackslash}m{#1}}

\begin{document}

\author{\normalsize
\mbox{\hspace{-2cm}
\begin{tabular}{ccc}
Eric Finster & Alex Rice & Jamie Vicary
\\
University of Birmingham & University of Cambridge & University of Cambridge
\\
e.l.finster@bham.ac.uk & alex.rice@cl.cam.ac.uk & jamie.vicary@cl.cam.ac.uk
\end{tabular}
\hspace{-2cm}}
}
\title{\bf A Type Theory for Strictly\\Associative Infinity Categories }

\maketitle

\begin{abstract}
  Many definitions of weak and strict $\infty$-categories have been proposed. In this paper we present a definition for $\infty$-categories with strict associators, but which is otherwise fully weak. Our approach is based on the existing type theory \Catt, whose models are known to correspond to weak $\infty$-categories. We add a definitional equality relation to this theory which identifies terms with the same associativity structure, yielding a new type theory  \(\Cattsa\), for strictly associative $\infty$-categories. We also provide a reduction relation which generates definitional equality, and show it is confluent and terminating, giving an algorithm for deciding equality of terms, and making typechecking decidable.

  Our key contribution, on which our reduction is based, is an operation on terms which we call \textit{insertion}. This has a direct geometrical interpretation, allowing a subterm to be inserted into the head of the term, flatting its syntactic structure. We describe this operation combinatorially in terms of pasting diagrams,  and also show can be characterized as a  pushout of contexts. This allows reasoning about insertion using just its universal property.

\end{abstract}

\section{Introduction}
\label{sec:intro}

It has long been known that not every fully weak \(\infty\)-category is equivalent to one which is fully strict\cite{simpson98_homot_types_stric}. While every bicategory (and by extension monoidal category) is equivalent to a strict 2-category (respectively strict monoidal category)\todo{citation?}, one can show that such an equivalence is not even possible at dimension \(3\), due to an idea involving the Eckmann-Hilton argument\cite{kock06_note_commut_doubl_semig_two}.

There are three main classes of coherence in a (globular) \(\infty\)-category: Unitors have the form
\[ f \cdot \id \Rightarrow f\]
associators have the form
\[ (f \cdot g) \cdot h \Rightarrow f \cdot (g \cdot h) \]
and interchangers have the form
\[ (a \cdot_n b) \cdot_m (c \cdot_n d) \Rightarrow (a \cdot_m c) \cdot_n (b \cdot_m d)\]
for compositions \(\cdot_n\) and \(\cdot_m\) of different codimensions (i.e. one could be a vertical composition and the other a horizontal composition). The Eckmann-Hilton argument implies that commutativity of scalars can be built from these 3 classes of coherences, and strictness of these coherences gives a strict commutativity of scalars, which in turn implies that not every weak \(\infty\)-category is equivalent to a strict one.

It is therefore clear that any definition of \(\infty\)-category in which all three of these classes is strict is not general enough to model all weak \(\infty\)-categories. However this does not mean that no structure of the category can be strictified, and this leads to the search of a semistrict definition of an \(\infty\)-category, one which is ``as strict as possible'' while still retaining the full expressivity of fully weak \(\infty\)-categories.

In dimension \(3\), two semistrict definitions have been given:
\begin{itemize}
\item Gray categories\cite{ae29f78ff3df4b86b920a204fff2eb68} are a definition where everything is strict apart from the interchange laws, which hold only up to isomorphism. It is proven that every tricategory is equivalent to a Gray category.
\item \citeauthor{joyalSimpson}\cite{joyalSimpson} introduce a definition for monoidal \(2\)-categories (i.e. \(3\)-categories with one object) which have weak units, but all other structure being strict, and show this is equivalent to the fully weak version.
\end{itemize}
These suggest two approaches to defining semistrict \(\infty\)-categories, one where unitors and associators are strict and interchangers are weak, as in Gray categories, and one where everything but units are strict, which is known as Simpson's conjecture.

Although many definitions of weak \(\infty\)-categories have been given (for example see \cite{leinster2001survey}), not many seem to be partially strict. Two of the authors, along with Reutter, have given a definition of a strictly unital \(\infty\)-category\cite{finster20_type_theor_stric_unital_categ} by putting a reduction relation on the type theory \(\Catt\)\cite{finster17_type_theor_defin_weak_ω}. Our work, running in parallel to this, uses the same technique but uses a new reduction relation which quotients out by associators instead of unitors.

\(\Catt\) is a type theory for which we can define the models to be \(\infty\)-categories. Recently it has been shown by \citeauthor{BenjaminCatt}\cite{BenjaminCatt}, that the models of this type theory are exactly Grothendieck-Maltsiniotis \(\infty\)-categories\cite{maltsiniotis2010grothendieck}. In contrast to the more abstract definitions, \(\Catt\) gives us a way to reason about the theory inductively and naturally allows computer assistance for defining and checking terms. By giving a complete reduction (one which is strongly normalising and confluent) we can create a quotient of this theory while retaining its computational properties such as decidable equality of terms and decidable typechecking.

To define our new reduction relation, we introduce an operation which we call ``insertion''. This operation effectively combines compound composites into a single composite. The simplest example of this is the standard 1-dimensional associativity equation. If we start with the compound composite \(f \cdot (g \cdot h)\), insertion will combine these into a single ternary composition \(f \cdot g \cdot h\). Similarly the term \((f \cdot g) \cdot h\) will be reduced to the same ternary composition. We call this operation insertion as when we view this as an operation on pasting diagrams, the pasting diagram for the inner composition is inserted into the pasting diagram for the outer composition. We give relatively weak conditions for this insertion to be admissible, and further show that it satisfies a certain universal property.

\subsection{Contributions and Future Work}
\label{sec:contrib}

We introduce a type theory \(\Cattsa\) for which the models form strictly associative \(\infty\)-categories. This is done by defining an operation on the structure of the type theory \(\Catt\) called insertion and creating an equality relation based on this. We further prove some properties of this insertion operation, and show that there is a reduction scheme which agrees with the equality and is complete, and hence give algorithms for deciding equality and typechecking.

Further we show that the pruning operation introduced for strict unit normalisation\cite{finster20_type_theor_stric_unital_categ} can be seen as an instance of our insertion operation, which suggests the two reduction schemes are compatible and can be combined to give a type theory which models strictly associative and unital \(\infty\)-categories. Proving that this reduction scheme is indeed complete is ongoing work.

Finally, we would eventually like to be able to find a connection between our work and associative n-categories\cite{dorn18_assoc}, which are a far more combinatoric definition of strictly associative and unital higher categories. These have the advantage that they can be viewed and manipulated graphically by contractions\cite{reutter2019high}, and have been implemented in the proof assistant ``homotopy.io''\cite{homotopyio-tool}. We eventually hope to find a strict enough type theoretic definition of higher categories such that associative n-category terms can be converted into the terms of this theory. This would have the two purposes: firstly it gives evidence that the structures manipulated by a tool such as ``homotopy.io'' really do correspond to conventional definitions of higher categories and secondly could make it easier to construct terms of \(\Catt\) and related theories.

\section{\textsf{Catt}}
\label{sec:catt}

We first define the syntax and typing rules of the base type theory \(\mathsf{Catt}\), which models weak \(\infty\)-categories. Here we give a full account of all the rules of the type theory, and carefully set up how we will add definitional equality into the theory, though do not give much motivation for many of the base roles, for which we encourage the reader to look at the original paper on \(\Catt\)\cite{finster17_type_theor_defin_weak_ω} or \citeauthor{benjamin:tel-03106197}'s PhD thesis\cite{benjamin:tel-03106197}.

\subsection{Syntax}
\label{sec:syntax}

Let \(V\) be an infinite set of variables. Then the syntax of
\emph{contexts}, \emph{substitutions}, \emph{types}, and \emph{terms}
(written \Ctx, \Sub, \Type, \Term) is given by:

\begin{center}
  \begin{tabular}{Sc Sc}
    {
    \begin{prooftree}
      \hypo{\phantom{\Term}} \infer1{\emptyset : \Ctx}
    \end{prooftree}
    }
    &
      {
      \begin{prooftree}
        \hypo{\Gamma : \Ctx} \hypo{A : \Type} \hypo{x \in V}
        \infer3{\Gamma, (x : A) : \Ctx}
      \end{prooftree}}
    \\
    {
    \begin{prooftree}
      \hypo{\phantom{\Term}} \infer1{\langle \rangle : \Sub}
    \end{prooftree}
    }
    & {
      \begin{prooftree}
        \hypo{\sigma : \Sub} \hypo{t : \Term} \hypo{x : V}
        \infer3{\langle \sigma , x \mapsto t \rangle : \Sub}
      \end{prooftree}
      }
    \\
    {
    \begin{prooftree}
      \hypo{\phantom{\Type}} \infer1{\star : \Type}
    \end{prooftree}
    }
    & {
      \begin{prooftree}
        \hypo{A : \Type} \hypo{s : \Term} \hypo{t : \Term}
        \infer3{\arr s A t : \Type}
      \end{prooftree}
      }
    \\
    {
    \begin{prooftree}
      \hypo{x \in V\vphantom{\Type}} \infer1{x : \Term}
    \end{prooftree}
    }
    & {
      \begin{prooftree}
        \hypo{\Gamma : \Ctx} \hypo{A : \Type} \hypo{\sigma : \Sub}
        \infer3{\Coh \Gamma A \sigma : \Term}
      \end{prooftree}
      }

  \end{tabular}
\end{center}
The symbol \(\equiv\) will represent syntactic equality.

We define the \emph{support} of a term, type, or substitution
in context \(\Delta\) as follows:
\begin{itemize}
\item Substitutions: \(\supp_\Delta(\langle \rangle) = \emptyset\),
  \(\supp_\Delta(\langle \sigma , t \rangle) = \supp_\Delta(\sigma)
  \cup \supp_\Delta(t)\).
\item Types: \(\supp_\Delta(\star) = \emptyset\), \(\supp(\arr s A t)
  = \supp_\Delta(s) \cup \supp_\Delta(t)\)
\item Terms:
  \begin{alignat*}{3}
    &\supp_\Delta(x) &&= \{x\} \cup \supp_\Delta(A) &\quad&\text{if \((x:A) \in \Delta\)}\\
    &\supp_\Delta(\Coh \Gamma A \sigma) &&= \supp_\Delta(\sigma)
  \end{alignat*}
\end{itemize}
We Write \(\supp(t)\) when the context is obvious.

A dimension function \(\dim\) on types is defined by:
\begin{alignat*}{2}
  &\dim(\star) &&= 0\\
  &\dim(\arr s A t) &&= \dim(A) + 1\\
  &\dim(\emptyset) &&= -1\\
  &\dim(\Gamma, (x : A)) &&= \max(\dim(\Gamma), \dim(A))
\end{alignat*}
And we define partial substitution functions for \(\sigma\, \tau :
\Sub, s\, t : \Term, A : \Type, x \in V\) as well as partial
composition function given by:
\begin{alignat*}{3}
  &\sub \star \sigma &&= \star\\
  &\sub {(\arr s A t)} \sigma &&= \arr {\sub s \sigma} {\sub A \sigma} {\sub t \sigma} &\quad&\text{if this is defined}\\
  &\sub {\Coh \Gamma A \tau} \sigma &&= \Coh \Gamma A {\tau \circ \sigma} &&\text{if \(\tau \circ \sigma\) is defined}\\
  &\sub x \sigma &&= t &&\text{if \(x \mapsto t \in \sigma\) and \(x\) only occurs once in \(\sigma\)}\\
  &\langle  \rangle \circ \sigma &&= \langle  \rangle\\
  &\langle \tau, x \mapsto t \rangle \circ \sigma &&= \langle \tau
  \circ \sigma, x \mapsto \sub t \sigma \rangle &&\text{if this is
    well defined}
\end{alignat*}

\subsection{Pasting Contexts and Globular Properties}
\label{sec:pds}

Contexts in \(\mathsf{Catt}\) correspond to finite computads (also known as finite polygraphs). Here we give rules to classify those contexts which form pasting diagrams.

\begin{center}
  \begin{tabular}{ScSc}
    {
    \begin{prooftree}
      \hypo{\vphantom{\Gamma \vdash_{pd} x : \star}}
      \infer1[$\star$]{x : \star \vdash_{pd} x : \star}
    \end{prooftree}
    }
    & {
      \begin{prooftree}
        \hypo{\Gamma \vdash_{pd} x : \star} \infer1[\checkmark]{\Gamma
          \vdash_{pd}}
      \end{prooftree}
      }
    \\
    {
    \begin{prooftree}
      \hypo{\Gamma \vdash_{pd} x : A} \hypo{y,f \in V \setminus
        \FV(\Gamma)} \infer2[$\Uparrow$]{\Gamma, (y : A), (f : \arr x
        A y) \vdash_{pd} f : \arr x A y}
    \end{prooftree}
    }
    & {
      \begin{prooftree}
        \hypo{\Gamma \vdash_{pd} f : \arr x A y}
        \infer1[$\Downarrow$]{\Gamma \vdash_{pd} y : A}
      \end{prooftree}
      }
  \end{tabular}
\end{center}

Then given a pasting scheme \(\Gamma\), the source
\(\partial^-(\Gamma)\) and target \(\delta^+(\Gamma)\) can be defined:
\begin{alignat*}{2}
  &\partial^-(x : \star) &&= \emptyset\\
  &\partial^-(\Gamma, y : A, f : \arr x A y) &&=
  \begin{cases*}
    \partial^-(\Gamma), y : A, f : \arr x A y&\text{if \(\dim(A) + 1 < \dim(\Gamma)\)}\\
    \partial^-(\Gamma)&\text{if \(\dim(A) + 1 = \dim(\Gamma)\)}\\
    \Gamma&\text{if \(\dim(A) + 1 > \dim(\Gamma)\)}
  \end{cases*}\\
  &\drop(\Gamma, x : A) &&= \Gamma\\
  &\partial^+(x : \star) &&= \emptyset\\
  &\partial^+(\Gamma, y : A, f : \arr x A y) &&=
  \begin{cases*}
    \partial^+(\Gamma), y : A, f : \arr x A y&\text{if \(\dim(A) + 1 < \dim(\Gamma)\)}\\
    \partial^+(\drop(\Gamma), y : A)&\text{if \(\dim(A) + 1 = \dim(\Gamma)\)}\\
    \drop(\Gamma), y : A&\text{if \(\dim(A) + 1 > \dim(\Gamma)\)}
  \end{cases*}
\end{alignat*}

More generally we have the notion of a \emph{globular} context, which is a context that does not contain any coherences in its syntax tree. These correspond to the computads which are just globular sets.

Given an ambient context \(\Gamma\), we can extend the dimension function to terms.

\begin{alignat*}{3}
  &\dim_\Gamma(x) &&= \dim(A) &\quad&\text{if \((x:A) \in \Gamma\)}\\
  &\dim_\Gamma(\Coh \Delta A \sigma) &&= \dim(A)\\
\end{alignat*}

Further for a term of dimension \(m\), and \(n \leq m\), we can define the \(n^{th}\) source and target of the term as follows.

\begin{alignat*}{3}
  &\delta_{\Gamma,n}^\epsilon(t) &&= t &\quad&\text{if \(n = \dim(x)\)}\\
  &\delta_{\Gamma,n}^\epsilon(x) &&= A^\epsilon_n &\quad&\text{if \((x: A) \in \Gamma \) and \(n < \dim(x)\)}\\
  &\delta_{\Gamma,n}^\epsilon(\Coh \Delta {A} \sigma) &&= \sub {A^\epsilon_n} \sigma &\quad&\text{if \(n < \dim(A)\)}\\
\end{alignat*}

As before, the subscript context will be omitted when it can be inferred.

We also introduce a way to deconstruct types:

\begin{definition}
  If \(A\) is a \(n\) dimensional type and \(m < n\) then define \(A^\epsilon_m\) as follows: First note that \(m < n\) so \(n\) is at least \(1\) and so \(A \equiv \arr s {A'} t\).
  \begin{itemize}
  \item If \(m = n - 1\) then \(A^-_m \equiv s\) and \(A^+_m \equiv t\).
  \item Otherwise \(A^\epsilon_m = {(A')}^\epsilon_m\).
  \end{itemize}
\end{definition}

Lastly we define unbiased compositions.

\begin{definition}
  Mutually inductively define the \emph{unbiased term} and \emph{unbiased type} of a pasting diagram \(\Gamma\), written \(UTm(\Gamma)\) and \(UTy(\Gamma)\).
  \begin{itemize}
  \item For a disk pasting diagram, the unbiased term is the top-dimensional variable, and the unbiased type is the type of this variable.
  \item For any other pasting diagram \(\Gamma\), the unbiased type is given by \(UTm(\partial^-(\Gamma)) \to UTm(\partial^+(\Gamma))\) and the unbiased term is given by \(\Coh \Gamma {UTy(\Gamma)} \id\)
  \end{itemize}

  A coherence \(\Coh \Delta A \sigma\) is \emph{unbiased} if \(A\) is the unbiased type over \(\Delta\).
\end{definition}

\subsection{Properties of Syntax}
\label{sec:props}

We will now introduce some auxiliary definitions on the syntax of \(\mathsf{Catt}\). To do we want to consider an arbitrary type theory based on the syntax. This must contain a notion of definitional equality, and typing judgments.

\begin{definition}
  A congruence \(\sim\) is a tuple of equivalence relations on contexts, terms, substitutions and types. By abuse of notation we call all 4 of these relations \(\sim\). We then must have the following rules:

  \begin{center}
  \begin{tabular}{Sc Sc}
    {
    \begin{prooftree}
      \hypo{\Gamma \sim \Gamma'} \hypo{A \sim A'}
      \infer2{\Gamma, (x : A) \sim \Gamma', (x : A')}
    \end{prooftree}
    }
    &
      {
      \begin{prooftree}
        \hypo{\sigma \sim \sigma'} \hypo{t \sim t'}
        \infer2{\langle \sigma, x \mapsto t \rangle \sim \langle \sigma', x \mapsto t' \rangle}
      \end{prooftree}
      }\\
    {
    \begin{prooftree}
      \hypo{A \sim A'} \hypo {s \sim s'} \hypo{t \sim t'}
      \infer3{\arr s A t \sim \arr {s'} {A'} {t'}}
    \end{prooftree}
    }&
       {
       \begin{prooftree}
         \hypo{\Gamma \sim \Gamma'} \hypo{A \sim A'} \hypo{\sigma \sim \sigma'}
         \infer3{\Coh \Gamma A \sigma \sim \Coh {\Gamma'} {A'} {\sigma'}}
       \end{prooftree}
       }
  \end{tabular}
\end{center}
\end{definition}

\begin{example}
  Syntactic equality, \(\equiv\), is a congruence.
\end{example}

We can now define the notion of a \(\mathsf{Catt}\)-based type theory.

\begin{definition}
  A \emph{\(\mathsf{Catt}\)-based type theory} consists of a congruence \(=\), and typing judgments:
  \begin{itemize}
  \item \(\Gamma \vdash\)
  \item \(\Gamma \vdash A\)
  \item \(\Gamma \vdash \sigma : \Delta\)
  \item \(\Gamma \vdash t : A\)
  \end{itemize}
  for contexts \(\Gamma\) and \(\Delta\), type \(A\), substitution \(\sigma\), and term \(t\). We further require that the following rules are satisfied:
  \begin{center}
  \begin{tabular}{Sc}
    {
    \begin{prooftree}
      \hypo{\Gamma \vdash} \hypo{\Gamma = \Gamma'}
      \infer2{\Gamma' \vdash}
    \end{prooftree}
    }\\
    {
    \begin{prooftree}
      \hypo{\Gamma \vdash A} \hypo{\Gamma = \Gamma'} \hypo{A = A'}
      \infer3{\Gamma' \vdash A}
    \end{prooftree}
    }\\
    {
    \begin{prooftree}
      \hypo{\Gamma \vdash \sigma : \Delta} \hypo{\Gamma = \Gamma'} \hypo{\sigma = \sigma'} \hypo{\Delta = \Delta'}
      \infer4{\Gamma' \vdash \sigma' : \Delta'}
    \end{prooftree}
    }\\
    {
    \begin{prooftree}
      \hypo{\Gamma \vdash t : A} \hypo{\Gamma = \Gamma'} \hypo{t = t'} \hypo{A = A'}
      \infer4{\Gamma' \vdash t' : A'}
    \end{prooftree}
    }
  \end{tabular}
\end{center}

\end{definition}

\subsection{Typing for \(\mathsf{Catt}\)}
\label{sec:typing}

To define the base type theory \(\mathsf{Catt}\), we need to define the definitional equality, and the typing judgments. The definitional equality will just be given by reflexivity. The typing rules for contexts, substitutions and types are given by:

\begin{center}
  \begin{tabular}{ScSc}
    {
    \begin{prooftree}
      \hypo{\vphantom{V \setminus \FV(\Gamma)}} \infer1{\emptyset
        \vdash}
    \end{prooftree}
    }
    & {
      \begin{prooftree}
        \hypo{\Gamma \vdash} \hypo{x \in V \setminus \FV(\Gamma)}
        \hypo{\Gamma \vdash A} \infer3{\Gamma, (x : A) \vdash}
      \end{prooftree}
      }
    \\
    {
    \begin{prooftree}
      \hypo{\Gamma \vdash} \infer1{\Gamma \vdash \star}
    \end{prooftree}
    }
    & {
      \begin{prooftree}
        \hypo{\Gamma \vdash A} \hypo{\Gamma \vdash s : A} \hypo{\Gamma
          \vdash t : A} \infer3{\Gamma \vdash \arr s A t}
      \end{prooftree}
      }
    \\
    {
    \begin{prooftree}
      \hypo{\Delta \vdash\vphantom{V \setminus \FV(\Gamma)}}
      \infer1{\Delta \vdash \langle \rangle : \emptyset}
    \end{prooftree}
    }
    & {
      \begin{prooftree}
        \hypo{\Delta \vdash \sigma : \Gamma} \hypo{\Gamma \vdash A}
        \hypo{\Delta \vdash t : \sub A \sigma} \hypo{x \in V \setminus
          \FV(\Gamma)} \infer4{\Delta \vdash \langle \sigma, x \mapsto
          t \rangle : \Gamma, (x : A)}
      \end{prooftree}
      }
  \end{tabular}
\end{center}

and then the rules for terms are given by:

\begin{center}
  \begin{tabular}{Sc}
    {
    \begin{prooftree}
      \hypo{\Gamma \vdash} \hypo{(x : A) \in \Gamma}
      \infer2[\((\text{var})\)]{\Gamma \vdash x : A}
    \end{prooftree}
    }\\
    {
    \begin{prooftree}
      \hypo{\Gamma \vdash_{pd}}
      \hypo{\Gamma \vdash \arr s A t}
      \hypo{\Delta \vdash \sigma : \Gamma}
      \hypo{\supp(s) = \partial^-(\Gamma)}
      \hypo{\supp(t) = \partial^+(\Gamma)}
      \infer5[\((\text{comp})\)]{\Delta \vdash \Coh \Gamma {(\arr s A t)} \sigma : \sub {(\arr s A t)} \sigma}
    \end{prooftree}
    }\\
    {
    \begin{prooftree}
      \hypo{\Gamma \vdash_{pd}} \hypo{\Gamma \vdash A} \hypo{\Delta
        \vdash \sigma : \Gamma} \hypo{\supp(A) = \Gamma}
      \infer4[\((\text{coh})\)]{\Delta \vdash \Coh \Gamma A \sigma : \sub A \sigma}
    \end{prooftree}
    }
  \end{tabular}
\end{center}

\begin{theorem}
  With a definitional equality given by syntactic equality, \(\mathsf{Catt}\) forms a \(\mathsf{Catt}\)-based type theory.
\end{theorem}

\section{Insertion}
\label{sec:insertion}

Now that we have defined the base type theory \(\Catt\), we are in a position to describe our insertion operation. We will first give an intuition to the action of this operation. Suppose we have the term \(a \cdot (b \cdot c)\), where \(a \cdot b\) represents (leaving out some implicit terms in the substitution):

\[ \Coh {(x : \star) (y : \star) (f : x \to y) (z : \star) (g : y \to z)} {x \to z} {f \mapsto a, g \mapsto b}\]

In this case we are composing over some pasting diagram, and then one of the (locally maximal) arguments to this pasting diagram is itself a composition over a pasting diagram, in this case the variable \(g\). Insertion then removes the variable \(g\) and everything in its support and ``inserts'' the inner pasting diagram in its place. This is demonstrated in \cref{fig:insertion} (where we have alpha renamed the inner pasting diagram to reduce confusion), where the final pasting diagram is simply giving the ternary composite, the unbiased version of the original composite.

To fully describe the new term, we also need two more constructions, as the coherence term constructor has three parts. Firstly we need a new substitution from the new pasting context to the ambient context, which assigns the arguments to the new composition. This can simply be done by mapping parts of the original pasting diagram (those in black in the figure) to their original arguments, and the parts of the inner pasting diagram (those in red in the figure) to the arguments of the inner composition. Secondly we need to construct the type part of the coherence (the \(A\) in \(\Coh \Gamma A \sigma\)). The most obvious way to do this is to transport the old type to the new pasting diagram, which will require a substitution from the old pasting diagram to the new inserted one.

\begin{figure}
  \centering
\[\begin{tikzcd}
        x & y & z && \textcolor{rgb,255:red,255;green,0;blue,0}{x'} & \textcolor{rgb,255:red,255;green,0;blue,0}{y'} & \textcolor{rgb,255:red,255;green,0;blue,0}{z'}
        \arrow["f", from=1-1, to=1-2]
        \arrow[""{name=0, anchor=center, inner sep=0}, "g", from=1-2, to=1-3]
        \arrow["{f'}", color={rgb,255:red,255;green,0;blue,0}, from=1-5, to=1-6]
        \arrow["{g'}", color={rgb,255:red,255;green,0;blue,0}, from=1-6, to=1-7]
        \arrow[curve={height=30pt}, shorten <=12pt, shorten >=12pt, Rightarrow, 2tail reversed, no head, from=0, to=1-6]
      \end{tikzcd}\]
    is sent to:
\[\begin{tikzcd}
        x & \textcolor{rgb,255:red,255;green,0;blue,00}{x'} & \textcolor{rgb,255:red,255;green,0;blue,0}{y'} & \textcolor{rgb,255:red,255;green,0;blue,0}{z'}
        \arrow["{f'}", color={rgb,255:red,255;green,0;blue,0}, from=1-2, to=1-3]
        \arrow["{g'}", color={rgb,255:red,255;green,0;blue,0}, from=1-3, to=1-4]
        \arrow["f", from=1-1, to=1-2]
\end{tikzcd}\]
  \caption{Insertion on \(a \cdot (b \cdot c)\)}
  \label{fig:insertion}
\end{figure}

Our insertion operation needs to trivialise associativity at all dimensions and all codimensions, and needs to consider composites of any shape so that the operation is natural enough to be confluent. Although the example we gave is relatively simple, describing the general operation at all dimensions can be tricky to reason about and so we will represent our pasting diagrams as labelled Batanin trees.

\begin{definition}
  Given a set of variables \(\L\), a \emph{\(\L\)-labelled Batanin trees}, \(T\) can be defined inductively as a tuple \((T_l,T_b)\) where \(T_l\) is a list of labels contained in \(\L\), \(T_b\) is a list of \(\L\) labelled Batanin trees, and further the length of \(T_l\) is one greater than the length of \(T_b\).
\end{definition}

\begin{definition}
  From every labelled Batanin tree, we can obtain a context, and further show that this context is a pasting diagram. Given a tree \(T\), we will denote the context \(\lfloor T \rfloor\).
\end{definition}

\begin{lemma}
  Given a context \(\Gamma\), with a proof that it is a pasting diagram, we can form a labelled Batanin tree \(\lceil \Gamma \rceil\). Further we have \(\lfloor \lceil \Gamma \rceil \rfloor = \Gamma\) and \(\lceil \lfloor T \rfloor \rceil = T\) for all pasting diagrams \(\Gamma\) and trees \(T\).
\end{lemma}
\begin{proof}
  See \cite{finster17_type_theor_defin_weak_ω}
\end{proof}

\begin{definition}
  The linear height of a Batanin tree is the height of the tree before any branches occur. It can be defined recursively as:
  \begin{equation*}
    \mathsf{lh}((T_l,T_b)) =
    \begin{cases}
      1 + \mathsf{lh}(t) &\text{if } T_b = [ t ],\\
      0 &\text{otherwise}
    \end{cases}
  \end{equation*}
\end{definition}

\begin{example}\label{ex:trees}
  A tree is generated from a pasting context by having the labels be the \(0\)-dimensional cells of the pasting diagram and then the tree between each 2 labels be the part of the pasting diagram suspended between those two labels. Consider these two examples, which we call \(\Delta\) and \(\Theta\):

  \begin{center}
    \begin{tabular}{|c|M{5cm}|M{5cm}|}
      \hline & Pasting diagram & Batanin tree\\ \hline
    \(\Delta\)& \raisebox{20pt}{\begin{tikzcd}[ampersand replacement = \&]
        x \&\& y \& z
        \arrow[from=1-3, to=1-4, "k"]
        \arrow[""{name=0, anchor=center, inner sep=0}, "g"{pos=0.7}, from=1-1, to=1-3]
        \arrow[""{name=1, anchor=center, inner sep=0}, "f"', curve={height=18pt}, from=1-1, to=1-3]
        \arrow[""{name=2, anchor=center, inner sep=0}, "h", curve={height=-18pt}, from=1-1, to=1-3]
        \arrow["\beta", shorten <=2pt, shorten >=2pt, Rightarrow, from=0, to=2]
        \arrow["\alpha", shorten <=2pt, shorten >=2pt, Rightarrow, from=1, to=0]
      \end{tikzcd}} &
                      {\begin{tikzpicture}[every node/.style={scale=0.8}]
                          \node [label={left:$x$},label={right:$z$},on grid](x01) {$\bullet$};
                          \node [label={left:$f$},label={right:$h$},above left=0.7 and 0.4 of x01, on grid] (x11) {$\bullet$};
                          \node [label={above:$k$},above right=0.7 and 0.4 of x01, on grid] (x12) {$\bullet$};
                          \node [label={above:$\alpha$},above left=0.7 and 0.25 of x11, on grid] (x21) {$\bullet$};
                          \node [label={above:$\beta$},above right=0.7 and 0.25 of x11, on grid] (x22) {$\bullet$};
                          \node at ($0.65*(x01) + 0.175*(x11) + 0.175*(x12)$) {$y$};
                          \node at ($0.65*(x11) + 0.175*(x21) + 0.175*(x22)$) {$g$};
                          \draw (x01.center) to (x11.center);
                          \draw (x01.center) to (x12.center);
                          \draw (x11.center) to (x21.center);
                          \draw (x11.center) to (x22.center);
                      \end{tikzpicture}}
    \\ \hline
\(\Theta\)&\raisebox{20pt}{\begin{tikzcd} [ampersand replacement=\&]
        {x'} \&\& {y'}
        \arrow[""{name=0, anchor=center, inner sep=0}, "{g'}"{pos=0.7}, from=1-1, to=1-3]
        \arrow[""{name=1, anchor=center, inner sep=0}, "{f'}"', curve={height=18pt}, from=1-1, to=1-3]
        \arrow[""{name=2, anchor=center, inner sep=0}, "{h'}", curve={height=-18pt}, from=1-1, to=1-3]
        \arrow["{\beta'}", shorten <=2pt, shorten >=2pt, Rightarrow, from=0, to=2]
        \arrow["{\alpha'}", shorten <=2pt, shorten >=2pt, Rightarrow, from=1, to=0]
      \end{tikzcd}}&
    {\begin{tikzpicture}[every node/.style={scale=0.8}]
                          \node [label={left:$x'$},label={right:$y'$},on grid](x01) {$\bullet$};
                          \node [label={left:$f'$},label={right:$h'$},above=0.7 of x01, on grid] (x11) {$\bullet$};
                          \node [label={above:$\alpha'$},above left=0.7 and 0.25 of x11, on grid] (x21) {$\bullet$};
                          \node [label={above:$\beta'$},above right=0.7 and 0.25 of x11, on grid] (x22) {$\bullet$};
                          \node at ($0.65*(x11) + 0.175*(x21) + 0.175*(x22)$) {$g'$};
                          \draw (x01.center) to (x11.center);
                          \draw (x11.center) to (x21.center);
                          \draw (x11.center) to (x22.center);
                      \end{tikzpicture}}\\ \hline
  \end{tabular}
\end{center}

  Here the first has a linear height of \(0\) and the second has a linear height of \(1\).
\end{example}

\begin{definition}
  Each label in a Batanin tree can be associated with a path (a non empty list of natural numbers) specifying how to reach that label from the root of the tree. A path in a Batanin tree is \emph{valid} if it is the path associated with one of the tree's labels.
\end{definition}

\begin{definition}
  Given a tree \(T\) and a locally maximal variable \(x\), we can define the branching path of this variable, which we think of as the path to the last branching in the tree before reaching this variable, as follows: If \(T_b\) is empty then we return the path \([0]\). Otherwise \(T_b\) is not empty and as \(x\) is locally maximal, it must occur in the \(n\)\textsuperscript{th} branch of \(T_b\) for some \(n\). If this branch is linear then we return \([n]\), otherwise we append \(n\) to the branching path of \(x\) in this branch.
\end{definition}

\begin{definition}[Insertion]
  Write \(A^{[i..j]}\) for the list consisting of the \(i\)\textsuperscript{th} element of \(A\) to the \(j\)\textsuperscript{th} element of \(A\) (zero indexing).
  Suppose \(S\) and \(T\) are Batanin trees, and \(p\) is a valid path in \(S\), and \(\mathsf{lh}(T) \geq \mathsf{length}(p) - 1\). Then we define the insertion of \(T\) into \(S\) along path \(p\), written \(\insertion S p T\), by recursion on \(p\).

  Suppose \(p = [n]\), then let \({(\insertion S p T)}_l\) be equal to \(S_l^{[0..(n-1)]} \doubleplus T_l \doubleplus S_l^{[(n+2)..]}\) and \({(\insertion S p T)}_b\) be equal to \(S_b^{[0..(n-1)]} \doubleplus T_b \doubleplus S_b^{[(n+1)..]}\).

  Suppose instead \(p\) has head \(n\) and non empty tail \(p'\). By the linearity condition we have that \(T_b = [T']\) for some \(t\) and let the \(n\)\textsuperscript{th} element of \(S_b\) be \(S'\). Then again let \({(\insertion S p T)}_l = S_l^{[0..(n-1)]} \doubleplus T_l \doubleplus S_l^{[n+2..]}\) and in this case let \({(\insertion S p T)}_b = S_b^{[0..(n-1)]} \doubleplus [ \insertion S' p' T' ] \doubleplus S_b^{[(n+1)..]}\).

  Finally for contexts \(\Delta\) and \(\Theta\), and some locally maximal variable \(x\) of \(\Delta\), let \(\insertion \Delta x \Theta = \lfloor \insertion {\lceil \Delta \rceil} p {\lceil \Theta \rceil} \rfloor\) where \(p\) is the branching path of \(x\).
\end{definition}

\begin{definition}
  The internal substitution from \(\iota_{\Delta,x,\Theta} : \Theta \to \insertion \Delta x \Theta\) is given my mapping each variable to itself.
\end{definition}

\begin{definition}
  The \emph{linear height} of a type \(A\) is the maximum \(n\) such that all terms contained in \(A\) with dimension less than or equal to \(A\) are variables.
\end{definition}
\begin{definition}
  Given a type \(A\) with \(\Delta \vdash A\) with linear height at least one less than the length of the branching path length of \(x\), we can define the external substitution \(\kappa_{\Delta,x,\Theta,A} : \Delta \to \insertion \Delta x \Theta\) by sending a variable \(y\) to \(y\) if variable \(y\) is in \(\insertion \Delta x \Theta\) and otherwise \(y\) is \(\delta^\epsilon_n(x)\) for some \(\epsilon \in \{+,-\}\) and \(n\), and then we can map \(y\) to \(\delta^\epsilon_n(\Coh {\Theta} A {\iota_{\Delta,x,\Theta}})\).
\end{definition}

\begin{figure}
  \centering
    \begin{tabular}{M{5cm}M{5cm}}
    \raisebox{20pt}{
\begin{tikzcd}[ampersand replacement=\&]
        \textcolor{rgb,255:red,255;green,0;blue,0}{x'} \&\& \textcolor{rgb,255:red,255;green,0;blue,0}{y'} \& z
        \arrow[from=1-3, to=1-4,"k"]
        \arrow[""{name=0, anchor=center, inner sep=0}, "{h'}"{pos=0.6}, color={rgb,255:red,255;green,0;blue,0}, curve={height=-12pt}, from=1-1, to=1-3]
        \arrow[""{name=1, anchor=center, inner sep=0}, "{f'}"', color={rgb,255:red,255;green,0;blue,0}, curve={height=30pt}, from=1-1, to=1-3]
        \arrow[""{name=2, anchor=center, inner sep=0}, "h", curve={height=-30pt}, from=1-1, to=1-3]
        \arrow[""{name=3, anchor=center, inner sep=0}, "{g'}"{pos=0.8}, color={rgb,255:red,255;green,0;blue,0}, curve={height=12pt}, from=1-1, to=1-3]
        \arrow["\beta", shorten <=2pt, shorten >=2pt, Rightarrow, from=0, to=2]
        \arrow["{\alpha'}", color={rgb,255:red,255;green,0;blue,0}, shorten <=2pt, shorten >=2pt, Rightarrow, from=1, to=3]
        \arrow["{\beta'}", color={rgb,255:red,255;green,0;blue,0}, shorten <=3pt, shorten >=3pt, Rightarrow, from=3, to=0]
\end{tikzcd}} &
                      {\begin{tikzpicture}[every node/.style={scale=0.8}]
                          \node [color={red},label={[red]left:$x'$},label={right:$z$},on grid](x01) {$\bullet$};
                          \node [color={red},label={[red]left:$f'$},label={right:$h$},above left=0.7 and 0.4 of x01, on grid] (x11) {$\bullet$};
                          \node [label={above:$k$},above right=0.7 and 0.4 of x01, on grid] (x12) {$\bullet$};
                          \node [color={red},label={[red]above:$\alpha'$},above left=0.7 and 0.5 of x11, on grid] (x21) {$\bullet$};
                          \node [label={above:$\beta$},above right=0.7 and 0.5 of x11, on grid] (x23) {$\bullet$};
                          \node [color={red},label={[red]above:$\beta'$},above=0.7 of x11, on grid](x22) {$\bullet$};
                          \node [red] at ($0.65*(x01) + 0.175*(x11) + 0.175*(x12)$) {$y'$};
                          \node [red] at ($0.5*(x11) + 0.2*(x21) + 0.3*(x22)$) {$g'$};
                          \node [red] at ($0.5*(x11) + 0.3*(x22) + 0.2*(x23)$) {$h'$};
                          \begin{pgfonlayer}{background}
                          \draw [red] (x01.center) to (x11.center);
                          \draw (x01.center) to (x12.center);
                          \draw [red] (x11.center) to (x21.center);
                          \draw [red] (x11.center) to (x22.center);
                          \draw (x11.center) to (x23.center);
                          \end{pgfonlayer}
                      \end{tikzpicture}}
    \end{tabular}
    \caption{The insertion of \(\Theta\) into \(\Delta\)}
  \label{fig:insert-ex}
\end{figure}
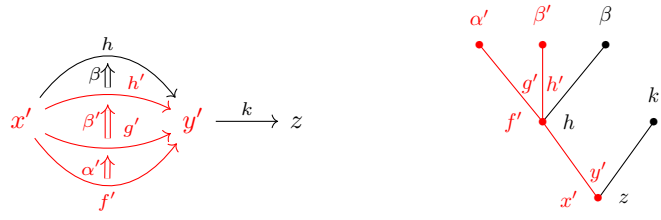

\begin{example}
  Suppose we have contexts \(\Delta\) and \(\Theta\) from \cref{ex:trees}. The variable \(\alpha\) has a branching height of \(1\) in \(\Delta\), and \(\Theta\) has a linear height of \(1\), so the insertion \(\insertion \Delta \alpha \Theta\) is well defined and is given in \cref{fig:insert-ex} as a pasting diagram and as a tree. The internal substitution \(\iota_{\Delta,\alpha,\Theta}\) maps \(\Theta\) into the part of \(\insertion \Delta \alpha \Theta\) marked in red by mapping each variable to itself.

  If we take \(A\) to be \(f \cdot k \to h \cdot k\), then \(\Delta \vdash A\) and \(A\) has linear height \(1\) and so we can define the external substitution \(\kappa_{\Delta,\alpha,\Theta,A}\) which has the following action:

  \begin{center}
  \begin{tabular}{ScSc}
    0-cells: & \(x \mapsto x'\quad y \mapsto y'\quad z \mapsto z\)\\
    1-cells: & \(f \mapsto f'\quad g \mapsto h'\quad h \mapsto h\quad k \mapsto k\)\\
    2-cells: & \(\alpha \mapsto \Coh \Theta A \iota\quad \beta \mapsto \beta\)
  \end{tabular}
\end{center}
\end{example}

\begin{definition}
  Given \(\sigma : \Delta \to \Gamma\), \(\tau : \Theta \to \Gamma\), and variable \(x\) in \(S\), we define the inserted substitution \(\insertion \sigma x \tau : \insertion \Delta {x} \Theta \to \Gamma\), as sending any variable \(y\) originating from \(S\) to \(\sigma(y)\) and any variable \(y\) originating from \(T\) to \(\tau(y)\).
\end{definition}

Subscripts from the external and internal substitutions may be dropped when there is no ambiguity.

\subsection{\(\Cattsa\)}
\label{sec:cattsa}

We define the type theory \(\Cattsa\) as having the same syntax as \(\Catt\). In addition to \(\Catt\) it will have an extra definitional equality and modify the typing rules for terms. These will be defined simultaneously, though this will be well founded as the typing for a term will only depend on equality for terms of lower dimension and equality of terms will only depend on typing of terms of the same dimension.

The definitional equality will be given by the smallest congruence such that for all well typed terms \(\Coh \Delta A \sigma\), locally maximal variable \(x\) of \(\Delta\) with \(\sigma(x)\) equal to an unbiased composite \(\Coh \Theta B \tau\) with the branching height of \(x\) is less than or equal to the linear height of \(\Theta\) we have:

\[\Coh \Delta A \sigma = \Coh {\insertion \Delta x \Theta} {\sub A {\kappa_{\Delta,x,\Theta}}} {\insertion \sigma x \tau} \]

Explicitly, the definitional equality will have the following structural rules:

\begin{center}
  \begin{tabular}{ScSc}
    {
    \begin{prooftree}
      \hypo{\vphantom{\Gamma}}
      \infer1{\emptyset = \emptyset}
    \end{prooftree}
    } &
        {
        \begin{prooftree}
          \hypo{\Gamma = \Gamma'} \hypo{A = A'}
          \infer2{\Gamma, x : A = \Gamma', x : A'}
        \end{prooftree}
        }\\
    {
    \begin{prooftree}
      \hypo{\vphantom{A = A'}}
      \infer1{\star = \star\vphantom{\arr s A t}}
    \end{prooftree}
    } & {
        \begin{prooftree}
          \hypo{A = A'} \hypo{s = s'} \hypo{t = t'}
          \infer3{\arr s A t = \arr {s'} {A'} {t'}}
        \end{prooftree}
        }\\
    {
    \begin{prooftree}
      \hypo{\vphantom{t'}}
      \infer1{\langle  \rangle = \langle  \rangle}
    \end{prooftree}
    } & {
        \begin{prooftree}
          \hypo{\sigma = \sigma'} \hypo{t = t'}
          \infer2{\langle \sigma, x \mapsto t \rangle = \langle \sigma', x \mapsto t' \rangle}
        \end{prooftree}
        }\\
    {
    \begin{prooftree}
      \hypo{t = t'}
      \infer1{t' = t}
    \end{prooftree}
    } & {
        \begin{prooftree}
          \hypo{t = t'} \hypo{t' = t''}
          \infer2{t = t''}
        \end{prooftree}
        }\\
    {
    \begin{prooftree}
      \hypo{\vphantom{A'}}
      \infer1{x = x\vphantom{\Coh Delta A \sigma}}
    \end{prooftree}
    } & {
        \begin{prooftree}
          \hypo{A = A'} \hypo{\sigma = \sigma'}
          \infer2{\Coh \Delta A \sigma = \Coh \Delta {A'} {\sigma'}}
        \end{prooftree}
        }
  \end{tabular}
\end{center}

\begin{prop}
  This definitional equality is a congruence.
\end{prop}

We then modify the typing rules for terms to the following:

\begin{center}
  \begin{tabular}{Sc}
    {
    \begin{prooftree}
      \hypo{\Gamma \vdash} \hypo{(x : A) \in \Gamma} \hypo{A = B}
      \infer3[\((\text{var}')\)]{\Gamma \vdash x : B}
    \end{prooftree}
    }\\
    {
    \begin{prooftree}
      \hypo{
        {
          \begin{matrix}
            \Gamma \vdash_{pd}&\Gamma \vdash \arr s A t&\Delta \vdash \sigma : \Gamma\\[3pt]
            \supp(s) = \partial^-(\Gamma)&\supp(t) = \partial^+(\Gamma)&B = \sub {(\arr s A t)} \sigma
          \end{matrix}
        }
      }
      \infer1[\((\text{comp}')\)]{\Delta \vdash \Coh \Gamma {(\arr s A t)} \sigma : B}
    \end{prooftree}
    }\\
    {
    \begin{prooftree}
      \hypo{\Gamma \vdash_{pd}}
      \hypo{\Gamma \vdash A}
      \hypo{\Delta \vdash \sigma : \Gamma}
      \hypo{\supp(A) = \Gamma}
      \hypo{B = \sub A \sigma}
      \infer5[\((\text{coh}')\)]{\Delta \vdash \Coh \Gamma A \sigma : B}
    \end{prooftree}
    }
  \end{tabular}
\end{center}

From these rules we can deduce that the following conversion rule is admissable:

\begin{prooftree*}
  \hypo{\Gamma \vdash t : A} \hypo{A = B}
  \infer2{\Gamma \vdash t : A'}
\end{prooftree*}

Next, we will prove that the constructions make some form of sense, in that the constructed substitutions are indeed well typed in \(\Cattsa\). We start with the following definition.

\begin{definition}
  A substitution \(\sigma\) from a globular set \(\Gamma\) to arbitrary context \(\Delta\) \emph{respects the globular conditions} if for all \(x : \arr s A t\) in \(\Gamma\), we have \(\delta^-(\sigma(x)) = \sigma(s)\) and \(\delta^+(\sigma(x)) = \sigma(t)\). Say a substitution is \emph{well-formed} if it respects the globular conditions and maps each variable of \(\Gamma\) to a term which is well typed in \(\Delta\).
\end{definition}

\begin{lemma}\label{lem:term-glob}
  If \(\Delta \vdash u : A\) then \(\delta^\epsilon_n(u) = A^\epsilon_n\) for all \(n < \dim(A)\).
\end{lemma}
\begin{proof}
  Proceed by induction on the proof of \(\Delta \vdash u : A\).
  \begin{itemize}
  \item If it is proved using a variable rule, then the result follows straight from definition of \(\delta\) and compatibility of \({(-)}^\epsilon_n\) and equality.
  \item If it is proved with the coherence or comp rule then \(u \equiv \Coh \Gamma B \sigma\) and \(A = \sub B \sigma\). Then \(\delta^\epsilon_n(u) \equiv \sub {B^\epsilon_n} \sigma \equiv {(\sub B \sigma)}^\epsilon_n = A^\epsilon_n\).
  \end{itemize}
\end{proof}


\begin{theorem}
  In \(\Cattsa\), if \(\Gamma\) is a globular set then \(\Delta \vdash \sigma : \Gamma\) if and only if \(\sigma\) is well-formed.
\end{theorem}
\begin{proof}
  Proceed by induction on the length of \(\sigma\). If \(\sigma\) is empty, then it is both well-formed and well-typed.

  Now consider \(\langle \sigma, x \mapsto u \rangle\) and first suppose \(\Delta \vdash \langle \sigma, x \mapsto u \rangle : \Gamma, x : A\). By inductive hypothesis, \(\sigma\) is well-formed. In particular each variable of \(\Gamma\) is sent by \(\sigma\) to a well-typed term of \(\Delta\). Further, by assumption \(\Delta \vdash u : \sub A \sigma\) and so \(x\) is also sent to a well-typed term. Now we need to show this substitution respects the globular conditions. The only case that does not immediately follow from inductive hypothesis is the one for \(x\). If \(A \equiv \star\) then there is nothing to do. Otherwise suppose \(A \equiv \arr s {A'} t\). Then \(\Delta \vdash u : \arr {\sub s \sigma} {\sub {A'} \sigma} {\sub t \sigma}\). Then by \cref{lem:term-glob} we have \(\delta^-(u) = \sub s \sigma\) and \(\delta^+(u) = \sub t \sigma\) as required.

  Now suppose \(\langle \sigma, x \mapsto u \rangle\) is well-formed. As a well-formed substitution maps all the variables in the source context, we must have \(x : A \in \Gamma\). Then we have that \(\Delta \vdash u : B\) for some type \(B\) by the substitution being well-formed. If \(B \equiv \star\) then \(A\) must also be \(\star\) so we are done. Otherwise, for all \(\epsilon\) and \(n\), we have that \(\delta^\epsilon_n(u) = B^\epsilon_n\) by \cref{lem:term-glob} and \(\delta^\epsilon_n(u) = \sub {A^\epsilon_n} \sigma\) by \(\sigma\) being well formed. Therefore \(\sub A \sigma = B\) and so \(\Delta \vdash u : \sub A \sigma\) by an application of the conversion rule.
\end{proof}

We now introduce a further equality operation, which effectively says that two terms are definitionally equal but also syntactically equal on any part of the term which has suitably high dimension. This will be useful later for proving confluence.

\begin{definition}
  Define the equality relation \(=_n\) on terms, types, and substitutions by the following rules.
  \begin{itemize}
  \item If \(s\) and \(t\) are terms with \(\dim(s) < n\) and \(\dim(t) < n\) then \(s =_n t\) if and only if \(s = t\).
  \item If \(s\) is a variable of dimension greater than or equal to \(n\) then \(s =_n t\) if and only if \(s \equiv t\).
  \item If \(s \equiv \Coh \Delta A \sigma\) with \(\dim(s) \geq n\) then \(s =_n t\) if and only if \(t \equiv \Coh {\Delta} {A'} {\sigma'}\) and \(A =_n A'\) and \(\sigma =_n \sigma'\).
  \item For types \(\star =_n \star\) and \(\arr s A t =_n \arr {s'} {A'} {t'}\) if and only if \(s =_n s'\), \(t =_n t'\), and \(A =_n A'\).
  \item For substitutions \(\langle t_0, \dots, t_k \rangle =_n \langle s_0, \dots, s_k \rangle\) if and only if \(t_i =_n s_i\) for all \(i\).
  \end{itemize}
\end{definition}

Next fix pasting contexts \(\Delta\) and \(\Theta\), a variable \(x\) in \(\Delta\), and a type \(A\) with \(\Theta \vdash A\), such that the linear height of \(A\) is at least one less than the branching height of \(x\).

\begin{prop}
  \(\iota : \Theta \to \insertion \Delta x \Theta\) is well-formed.
\end{prop}
\begin{proof}
  As all the variables in \(\Theta\) have the same type in the inserted context, this is clearly a well-formed substitution.
\end{proof}


\begin{prop}
  The substitution \(\kappa : \Delta \to \insertion \Delta x \Theta\) is well-formed.
\end{prop}
\begin{proof}
  For this to be a well-formed substitution, the only interesting case is where \(y\) is not in the support of \(x\), yet one or both of it's source and target variables are. Suppose variable \(z\) is the target of \(y\) and the \(n\)\textsuperscript{th} source of \(x\) in \(\Delta\). An inspection of the algorithm shows that the target of \(y\) in \(\insertion \Delta x \Theta\) is the unique variable in the \(n\)\textsuperscript{th} source of \(\Theta\), which we call \(z'\). However now we have that \(\kappa(z) \equiv \delta^-_n(\Coh \Theta A \iota) \equiv z' \equiv \delta^-(y)\) and so \(\kappa(\delta_\Delta^-(y)) \equiv \delta_{\insertion \Delta x \Theta}^-(\kappa(y))\). The case where \(z\) was the source of \(y\) follows similarly.
\end{proof}

Now further assume that we have well-formed \(\sigma : \Delta \to \Gamma\) and \(\tau : \Theta \to \Gamma\) with \(\sigma(x) = \Coh \Theta A \tau\).

\begin{prop}\label{prop:sub-comm}
  \(\tau = (\insertion \sigma x \tau) \circ \iota\) and \(\sigma = (\insertion \sigma x \tau) \circ \kappa\). 
\end{prop}
\begin{proof}
  The first equality follows directly from the definition of \(\insertion \sigma x \tau\). To prove the second, it suffices to prove the cases for \(y\) in the support of \(x\), as the rest hold definitionally. If \(y = \delta^{\epsilon}_n(x)\), then \(\kappa(y) = \delta^{\epsilon}_n(\Coh \Theta A {\iota})\). Further we know that, as \(\sigma(x) = \Coh \Theta A \tau\) and as all the substitutions are well-formed, \(\sub y \sigma\) is equal to \(\delta^\epsilon_n(\Coh \Theta A \tau)\). Therefore:
  \begin{align*}
    \sub y \sigma &= \delta_n^\epsilon(\Coh \Theta A \tau)\\
                  &= \sub {\delta_n^\epsilon(\Coh \Theta A \id)} {\tau}\\
                  &= \sub {\delta_n^\epsilon(\Coh \Theta A \id)} {(\insertion \sigma x \tau) \circ \iota}\\
                  &= \sub {(\sub {\delta_n^\epsilon(\Coh \Theta A \id)} {\iota})} {\insertion \sigma x \tau}\\
                  &= \sub {(\delta_n^\epsilon(\Coh \Theta A {\iota}))} {\insertion \sigma x \tau}\\
                  &= \sub {(\sub y {\kappa})} {\insertion \sigma x \tau} \\
                  &= \sub y {(\insertion \sigma x \tau) \circ \kappa}
  \end{align*}
  and so we get the required syntactic equality.
\end{proof}

\begin{prop}
  \(\insertion \sigma x \tau\) is well-formed.
\end{prop}
\begin{proof}
  To check that \(\insertion \sigma x \tau\) is well-formed, we again look at the one problematic case where \(y\) originates from \(\Delta\) yet its (without loss of generality) target \(z\) originates from \(\Theta\). Then we need to show that \(\delta^+(\sigma(y)) = \tau(z)\). But then as \(\kappa\) is a well-formed substitution we have that there is \(z' = \delta_\Delta^+(y)\) with \(\kappa(z') = z\). As \(\sigma\) is well-formed we have that \(\delta^+(\sigma(y)) = \sigma(z') = \sub {(\kappa(z'))} {\insertion \sigma x \tau} = (\insertion \sigma x \tau)(z) = \tau(z)\) as required.
\end{proof}

\subsection{Universal Property of Insertion}
\label{sec:uprop}

We first must introduce some theory about disk contexts:

\begin{definition}
  Recursively define the disk contexts \(D_n\) as \(D_0\) being the context \(d_0^- : \star\) and \(D_{n+1}\) being the context \(D_n, d_n^+ : d_{n-1}^- \to d_{n-1}^+, d_{n+1}^- : d_n^- \to d_n^+\).
\end{definition}

Note that each disk context has a single maximal dimension cell, and any substitution out of a disk context is fully determined by where this cell is sent.

\begin{definition}
  A substitution \(D_n \to \Gamma\) is fully determined by its action on \(d_n^-\). If \(\Gamma \vdash t\) and \(t\) is \(n\)-dimensional, then let \(\overline t : D_n \to \Gamma\) be the substitution that sends \(d_m^\epsilon\) to \(\delta_m^\epsilon(t)\) (and so in particular sends \(d_n^-\) to \(t\)).
\end{definition}

Given pasting contexts \(\Delta\) and \(\Theta\), a variable \(x\) in \(\Delta\) and a type \(A\) with \(\Theta \vdash A\), we can classify insertion by a universal property.

\begin{definition}
  When \(\dim(A) = \dim(x) = n\), we say that an insertion of \((\Theta,A)\) into \(\Delta\) via \(x\) is the following pushout:
  \[\begin{tikzcd}
      {D_n} & \Delta \\
      \Theta & I
      \arrow["{\overline {\Coh \Theta A \id}}"', from=1-1, to=2-1]
      \arrow["\overline x"', from=1-1, to=1-2]
      \arrow[from=2-1, to=2-2]
      \arrow[from=1-2, to=2-2]
      \arrow["\lrcorner"{anchor=center, pos=0.125, rotate=180}, draw=none, from=2-2, to=1-1]
    \end{tikzcd}\]
\end{definition}

Note that as a pushout, this notion is unique up to isomorphism. If we further require that the insertion is a pasting diagram, then (if it exists) this uniquely specifies this colimit.

\begin{theorem}
  When \(\dim(A) = \dim(x) = n\), \(x\) is locally maximal, and \(A\) is linear up to at least one less than the length of the branching path of \(x\), then the following is a pushout in \(\Cattsa\):
  \[\begin{tikzcd}
      {D_n} & \Delta \\
      \Theta & {\insertion \Delta x \Theta}
      \arrow["{\overline {\Coh \Theta A \id}}"', from=1-1, to=2-1]
      \arrow["\overline x"', from=1-1, to=1-2]
      \arrow["{\iota}", from=2-1, to=2-2]
      \arrow["{\kappa}"', from=1-2, to=2-2]
      \arrow["\lrcorner"{anchor=center, pos=0.125, rotate=180}, draw=none, from=2-2, to=1-1]
    \end{tikzcd}\]
  and hence \(\insertion \Delta x \Theta\) is an insertion of \((\Theta,A)\) into \(\Delta\) along \(x\). Further it is a pasting diagram and so is the unique pasting diagram with this property.
\end{theorem}
\begin{proof}
  We have already checked that all these contexts and substitutions are well typed. To check commutativity of the square it suffices to check that the maximal variable is sent to the same term by each substitution. Along the top it is sent to \(\sub x {\kappa}\), which is defined to be \(\Coh{\Theta} A {\iota} = \sub {(\Coh \Theta A \id)} {\iota}\), which is the term it is sent to along the bottom.

  Now suppose that we have \(\Gamma\), \(\sigma\), and \(\tau\) such that the following commutes:
  \[
    \begin{tikzcd}
      {D_n} & \Delta \\
      \Theta & {\insertion \Delta x \Theta} \\
      && \Gamma
      \arrow["{\Coh \Theta A \id}"', from=1-1, to=2-1]
      \arrow["x"', from=1-1, to=1-2]
      \arrow["{\iota}", from=2-1, to=2-2]
      \arrow["{\kappa}"', from=1-2, to=2-2]
      \arrow["\sigma", curve={height=-18pt}, from=1-2, to=3-3]
      \arrow["\tau"', curve={height=18pt}, from=2-1, to=3-3]
    \end{tikzcd}
  \]

  Then we need to generate a substitution \(\insertion \Delta x \Theta \to \Gamma\), for which we can give \(\insertion \sigma x \tau\). Any variable \(l\) that originates from \(\Delta\) must be sent to \(\sigma(l)\) and any variable \(l\) that originates from \(\Theta\) must be sent to \(\tau(l)\). This ensures uniqueness. The commutativity conditions are given in \cref{prop:sub-comm}.
\end{proof}

To finish this section we note that the conditions we have required for insertion to work are much weaker than what we have used for the definition of \(\Cattsa\). Firstly, for the equality in \(\Cattsa\), we require that the inner term be fully unbiased instead of its type being linear up to a certain height, and further the insertion operation does not depend on the inner term being a composition instead of a coherence, which we explore in the next example.

\begin{example}
  Suppose we have pasting diagram \(\Delta\) with locally maximal \(x\), and \(\sigma : \Delta \to \Gamma\) with \(\sigma(x) \equiv \id(t)\). Then the pasting diagram of \(\id(t)\) is \(D_n\) for some \(n\). As \(D_n\) is fully linear, it satisfies the conditions for insertion, and so we can form the insertion \(\insertion \Delta x {D_n}\) and from a term \(\Coh \Delta A \sigma\), we could form the term \(\Coh {\insertion \Delta x {D_n}} {\sub A \kappa} {\insertion \sigma x {\overline t}}\). This term has the same form\footnote{The two operations are not identical, as inserting an identity takes as many arguments from the inner term as possible, and pruning discards these arguments, taking them from the outer term instead. However, the normal forms of these two operations (with appropriate structural rules) should be identical.} as the result of applying the pruning operation from \cite{finster20_type_theor_stric_unital_categ} to the original term.
\end{example}

That pruning can be realised as a case of insertion has two exciting consequences. Firstly, this means that we can also realise pruning as a universal property, which was not previously known and allows some proofs to be written in a neat way. It also suggests that the two operations should be compatible and it is possible that they can be combined to get a strictly associative and unital type theory. Proving confluence of such a system is ongoing work.

\subsection{Properties of equality}
\label{sec:prop-eq}

We can now prove some facts about the equality of \(\Cattsa\).

\begin{lemma}\label{lem:inserted-sub}
  Suppose there exists an insertion of \((\Theta,A)\) into \(\Delta\) and we have substitutions \(\rho : \Delta \to \Gamma\), \(\tau : \Theta \to \Gamma\) and \(\sigma : \Gamma \to \Gamma'\) with \(\rho(x) = \Coh \Theta A \tau\). Then \(\sigma \circ (\insertion \rho x \tau) = \insertion {(\sigma \circ \rho)} x {(\sigma \circ \tau)}\).
\end{lemma}
\begin{proof}
  Take \(y \in \Var(\insertion \Delta x \Theta)\). If \(y \in \Var(\Theta)\) then:
  \begin{align*}
    (\sigma \circ (\insertion \rho x \tau))(y) &= \sub {((\insertion \rho x \tau)(y))} {\sigma}\\
                                               &= \sub {(\tau(y))} {\sigma}\\
                                               &= (\sigma \circ \tau)(y)\\
                                               &= (\insertion {(\sigma \circ \rho)} x {\sigma \circ \tau})(y)
  \end{align*}
  If \(y \in \Var(\Delta)\) then:
  \begin{align*}
    (\sigma \circ (\insertion \rho x \tau))(y) &= \sub {((\insertion \rho x \tau)(y))} {\sigma}\\
                                               &= \sub {(\rho(y))} {\sigma}\\
                                               &= (\sigma \circ \rho)(y)\\
                                               &= (\insertion {(\sigma \circ \rho)} x {\sigma \circ \tau})(y)
  \end{align*}
  As these are the only two cases, this completes the proof.
\end{proof}

\begin{lemma}\label{lem:sub-equality}
  For all well typed \(\sigma = \sigma'\) we have:
  \begin{itemize}
  \item \(\sub t \sigma = \sub {t'} {\sigma'}\) for all well typed terms \(t = t'\)
  \item \(\sub A \sigma = \sub {A'} {\sigma'}\) for all well typed types \(A = A'\)
  \item \(\sigma \circ \tau = \sigma \circ \tau'\) for all well typed substitutions \(\tau = \tau'\).
  \end{itemize}
\end{lemma}
\begin{proof}
  Prove by mutual induction on the equality of types, terms, and substitutions.
  \begin{itemize}
  \item Take equality \(A = A'\). If this is an instance of \(\star = \star\) then \(\sub \star \rho \equiv \star\) for any \(\rho\), so the result is trivial. If it is instead an instance of \(\arr s B t = \arr {s'} {B'} {t'}\) from \(s = s'\), \(t = t'\), and \(B = B'\). Then by inductive hypothesis we get:
    \begin{align*}
      \sub A \sigma &= \sub {(\arr s B t)} \sigma\\
                    &= \arr {\sub s \sigma} {\sub B \sigma} {\sub t \sigma}\\
                    &= \arr {\sub {s'} {\sigma'}} {\sub {B'} {\sigma'}} {\sub {t'} {\sigma'}}\\
                    &= \sub {\arr {s'} {B'} {t'}} {\sigma'}\\
                    &= \sub {A'} {\sigma'}
    \end{align*}
  \item Take equality \(t = t'\). If this is an instance of \(x = x\) then we need that \(\sigma(x) = \sigma'(x)\), which follows from \(\sigma = \sigma'\). If it is an instance of transitivity or symmetry, then we can simply use inductive hypothesis and transitivity or symmetry on the resulting terms. For the structural coherence case we have \(\Coh \Delta B \rho = \Coh \Delta {B'} {\rho'}\) with \(B = B'\) and \(\rho = \rho'\). Then \(\sub t \sigma \equiv \Coh \Delta B {\sigma \circ \rho}\) and \(\sub {t'} {\sigma'} \equiv \Coh \Delta {B'} {\sigma' \circ \rho'}\). Using inductive hypothesis we get that \(\sigma \circ \rho = \sigma' \circ \rho'\) and so we can apply the structural coherence rule to get the required equality.

    This leaves the insertion case. Suppose we have \(\Coh \Delta B \rho\) with \(\rho(x) = \Coh \Theta C \tau\) and the equality is \(\Coh \Delta B \rho = \Coh {\insertion \Delta x \Theta} {\sub A {\kappa_{\Delta,x,\Theta}}} {\insertion \rho x \tau} \). When we apply the substitutions we get that \(\sub t \sigma \equiv \Coh \Delta B {\sigma \circ \rho}\) and \(\sub {t'} {\sigma'} \equiv \Coh {\insertion \Delta x \Theta} {\sub B {\kappa_{\Delta,x,\Theta}}} {\sigma' \circ (\insertion \rho x \tau)}\). As \((\sigma \circ \rho) (x) \equiv \Coh \Theta C {\sigma \circ \tau}\), we can certainly apply an insertion to \(\sub t \sigma\) to get \(\Coh {\insertion \Delta x \Theta} {\sub B {\kappa_{\Delta,x,\Theta}}} {\insertion {(\sigma \circ \rho)} x {(\sigma \circ \tau)}}\).

    It remains to show modulo structural rules that \(\sigma' \circ (\insertion \rho x \tau) = \insertion {(\sigma \circ \rho)} x {(\sigma \circ \tau)}\). By inductive hypothesis we have that \(\sigma' \circ (\insertion \rho x \tau) = \sigma \circ (\insertion \rho x \tau)\), and the remaining equality is given by \cref{lem:inserted-sub}.
  \item The case for substitutions follows easily from inductive hypothesis.
  \end{itemize}
\end{proof}

\begin{lemma}\label{lem:support-equality}
  The support operation is invariant under equality.
\end{lemma}

\begin{prop}
  \(\Cattsa\) is a \(\Catt\)-based type theory.
\end{prop}

\begin{proof}
  We proceed by mutually proving the following statements by induction on the equality proofs.
  \begin{itemize}
  \item For any \(\Gamma = \Gamma'\) we have \(\Gamma \vdash\) if and only if  \(\Gamma' \vdash\).
  \item For any \(A = A'\) we have that for any \(\Gamma = \Gamma'\) that \(\Gamma \vdash A\) if and only if \(\Gamma' \vdash A'\).
  \item For any \(\sigma = \sigma'\), for all \(\Gamma = \Gamma'\) and \(\Delta = \Delta'\), \(\Gamma \vdash \sigma : \Delta\) if and only if \(\Gamma' \vdash \sigma' : \Delta'\).
  \item For any \(t = t'\), for all \(\Gamma = \Gamma'\) and \(A = A'\), \(\Gamma \vdash t : A\) if and only if \(\Gamma' \vdash t' : A'\).
  \end{itemize}
  Now we split into cases:
  \begin{itemize}
  \item Suppose we have \(\Gamma = \Gamma'\) and want to prove that \(\Gamma \vdash\) if and only if \(\Gamma' \vdash\). Case splitting on the equality:
    \begin{itemize}
    \item If \(\Gamma = \Gamma'\) is an instance of \(\emptyset = \emptyset\) then the result is trivial.
    \item Otherwise we have \(\Delta, x : A = \Delta, x : A'\) with \(\Delta = \Delta'\) and \(A = A'\). If \(\Delta, x : A \vdash\) then \(\Delta \vdash\) and \(\Delta \vdash A\). By inductive hypothesis, we have that \(\Delta' \vdash\) and \(\Delta' \vdash A'\) and so \(\Delta', x : A' \vdash\). The other direction is symmetrical.
    \end{itemize}
  \item Suppose we have \(A = A'\). We want to show for all \(\Gamma = \Gamma'\) that \(\Gamma \vdash A\) if and only if \(\Gamma' \vdash A'\). Case splitting on \(A = A'\).
    \begin{itemize}
    \item If it is a case of \(\star = \star\) then the result is trivially true.
    \item Otherwise the equality is an instance of \(\arr s A t = \arr {s'} {A'} {t'}\) with \(s = s'\), \(A = A'\), \(t = t'\). As above, the result follows easily from inductive hypothesis.
    \end{itemize}
  \item Suppose we have \(\sigma = \sigma'\). We want to show for all \(\Gamma = \Gamma'\) and \(\Delta = \Delta'\) that \(\Gamma \vdash \sigma : \Delta\) if and only if \(\Gamma' \vdash \sigma' : \Delta'\): Cases splitting on \(\sigma = \sigma'\):
    \begin{itemize}
    \item If the equality is an instance of \(\langle  \rangle = \langle  \rangle\) then again the result is trivial.
    \item If the equality is an instance of \(\langle \tau, x \mapsto t \rangle = \langle \tau, x \mapsto t' \rangle\) with \(\tau = \tau'\) and \(t = t'\) and then suppose \(\Gamma \vdash \sigma : \Delta\). Then it must be that \(\Delta \equiv \Theta , x : A\) and so as \(\Delta = \Delta'\), \(\Delta' \equiv \Theta', x : A'\) with \(\Theta = \Theta'\) and \(A = A'\). We further have that \(\Gamma \vdash \tau : \Theta\) and \(\Gamma \vdash t : \sub A \tau\). By \cref{lem:sub-equality}, \(\sub A \tau = \sub {A'} {\tau'}\) and so by inductive hypothesis we have that \(\Gamma' \vdash \tau' : \Theta'\) and \(\Gamma' \vdash t' : \sub {A'} {\tau'}\) and so \(\Gamma' \vdash \sigma' : \Delta'\) as required. The other direction again follows the same argument.
    \end{itemize}
  \item Suppose we have \(t = t'\). We want to show for all \(\Gamma = \Gamma'\) and \(A = A'\) that \(\Gamma \vdash t : A\) if and only if \(\Gamma' \vdash t' : A'\). Case split on \(t = t'\):
    \begin{itemize}
    \item If it is a case of \(x = x\) then suppose we have \(\Gamma \vdash x : A\). Then this typing rule must be derived from (var') so \(x : B \in \Gamma\) and \(A = B\). Then \(x : B' \in \Gamma\) with \(B = B'\) as \(\Gamma = \Gamma'\) and so as equality is a congruence we have \(A' = A = B = B'\) and so we can apply the (var') rule to get \(\Gamma' \vdash x : A'\). The other direction is symmetrical.
    \item If it is a case of symmetry, then the result follows directly from inductive hypothesis.
    \item If it is a case of transitivity, then the result also follows directly from inductive hypothesis.
    \item If it is from the structural rule on coherences we have \(\Coh \Delta B \sigma = \Coh \Delta {B'} {\sigma'}\) with \(B = B'\) and \(\sigma = \sigma'\). Suppose we have \(\Gamma \vdash \Coh \Delta B \sigma : A\) by the composition rule or coherence rule. In this case we have \(A = \sub B \sigma\), \(\Delta \vdash B\), and \(\Gamma \vdash \sigma : \Delta\). Then by inductive hypothesis we have \(\Delta \vdash B'\) and \(\Gamma' \vdash \sigma' : \Delta\). Further we have \(A' = A = \sub B \sigma = \sub {B'} {\sigma'}\) by equality being a congruence and \cref{lem:sub-equality}. Then, using \cref{lem:support-equality} and a coherence or composition rule, we get that \(\Gamma' \vdash t' : A'\). The other direction follows similarly.
    \item The last case is where the equality is an instance of insertion. Here we have \(t \equiv \Coh \Delta B \sigma\) with \(\sigma(x) = \Coh \Theta C \tau\) and \(t' \equiv \Coh {\insertion \Delta x \Theta} {\sub B {\kappa}} {\insertion \sigma x \tau}\).

      For the reverse direction we assume that \(\Gamma \vdash t' : A'\). By assumption of insertion we have that \(\Gamma \vdash t : D\) for some type \(D\). Then by inspection of rules we must have that \(D = \sub B \sigma\). Similarly we have that \(A' = \sub {\sub B \kappa} {\insertion \sigma x \tau} = \sub B {(\insertion \sigma x \tau) \circ \kappa} = \sub B \sigma\) and so \(D = A'\).

      Conversely assume that \(\Gamma \vdash t : A\). Using inductive hypothesis we get that \(B\) is well typed in \(\Delta\) and \(\sigma\) is well typed from \(\Delta\) to \(\Gamma'\). Then by results in the previous section we have that \(\insertion \Delta x \Theta\) is a pasting context, \(\kappa\) is well typed and so \(\sub B \kappa\) is well typed, and \(\insertion \sigma x \tau\) is a valid substitution from \(\insertion \Delta x \Theta\) to \(\Gamma'\). Therefore \(t'\) can be typed with any type equal to \(\sub {\sub B \kappa} {\insertion \sigma x \tau} = \sub B \sigma\). Now suppose \(A' = A\). Similarly to before we have \(A = \sub B \sigma\) so by \(=\) being a congruence we have \(A' = \sub B \sigma\) and so \(\Gamma' \vdash t' : A'\).
    \end{itemize}
  \end{itemize}
  This completes all cases and is sufficient to show what was required.
\end{proof}

\section{Reduction}
\label{sec:reduction}

So far we have given a type theory which models strictly associative infinity categories. Unfortunately, our type theory currently has no way to decide equality and by extension perform type checking. To solve this we introduce a reduction relation, show that it agrees with the definitional equality, and show that it is complete by showing it has strong normalisation and local confluence.

\begin{definition}[Reduction on types]
  The following rules are admissible for reducing types.
  \begin{itemize}
  \item If \(s \rto s'\):
    \[
      \arr s A t \rto \arr {s'} A t
    \]
  \item If \(t \rto t'\):
    \[
      \arr s A t \rto \arr s A {t'}
    \]
  \item If \(A \rto A'\):
    \[
      \arr s A t \rto \arr s {A'} t
    \]
\end{itemize}
\end{definition}

\begin{definition}[Reduction on substitutions]
  The following rules can be used to reduce substitutions:
  \begin{itemize}
  \item If \(\sigma \rto \sigma'\):
    \(\langle \sigma, t \rangle \rto \langle \sigma', t \rangle\)
  \item If \(t \rto t'\):
    \(\langle \sigma, t \rangle \rto \langle \sigma, t' \rangle\)
  \end{itemize}
  This effectively allows any term in the substitution to be reduced.
\end{definition}

\begin{definition}[Reduction on terms]
  The following rules can be used to reduce terms:
  \begin{itemize}
  \item (Insertion): If \(x\) is locally maximal in \(\Delta\) and \(\sigma(x)\) is the unbiased composite \(\Coh \Theta B \tau\) with the branching height of \(x\) being less than or equal to the linear height of \(\Theta\) then:
    \[\Coh \Delta A \sigma \rto \Coh {\insertion \Delta x \Theta} {\sub A {\kappa_{\Delta,x,\Theta}}} {\insertion \sigma x \tau} \]
  \item (Cell reduction): The type of a coherence can be reduced, so if \(A \rto A'\) then:
    \[\Coh {\Delta} A \sigma \rto \Coh {\Delta} {A'} \sigma\]
  \item (Argument reduction): The arguments to a coherence can be reduced, so if \(\sigma \rto \sigma'\) then:
    \[\Coh \Delta A \sigma \rto \Coh \Delta A {\sigma'}\]
  \end{itemize}
\end{definition}

\begin{lemma}\label{lem:red-sub}
  For all terms \(t \rto t'\), types \(A \rto A'\), and substitutions \(\sigma \rto \sigma'\), and another substitution \(\tau\), then \(\sub t \tau \rto \sub {t'} \tau\), \(\sub A \tau \rto \sub {A'} \tau\) and \(\tau \circ \sigma \rto \tau \circ \sigma'\).
\end{lemma}
\begin{proof}
  Proceed by mutual induction on \(t\), \(A\), and \(\sigma\).
  \begin{itemize}
  \item For substitutions \(\sigma = \langle t_0, \dots, t_n \rangle\), we have \(\sigma \rto \sigma'\) and therefore there must be some \(i\) with \(t_i \rto t_i'\) and \(\sigma' = \langle t_0,\dots,t_i',\dots,t_n \rangle\). Further, \(\tau \circ \sigma = \langle \sub {t_0} \tau, \dots, \sub {t_n} \tau \rangle\) and \(\tau \circ \sigma' = \langle \sub {t_0} \tau, \dots, \sub {t_i'} \tau, \dots, \sub {t_n} \tau \rangle\) and so by inductive hypothesis \(\sub {t_i} \tau \rto \sub {t_i'} \tau\) and so \(\tau \circ \sigma \rto \tau \circ \sigma'\).
  \item For types we get a proof similar to above where we apply the inductive hypothesis to some subterm.
  \item For terms \(t\), we note that \(t\) cannot be a variable as variables do not reduce. Therefore \(t \equiv \Coh \Delta A \sigma\) and \(\sub t \tau \equiv \Coh \Delta A {\tau \circ \sigma}\). Now we case split on the possible reductions of \(t\).
    \begin{itemize}
    \item Cell reduction: The type \(A\) is unaffected by substitution so this remains unchanged.
    \item Argument reduction: Follows from inductive hypothesis.
    \item Insertion: Suppose \(t \rto t'\) is an insertion along \(x\) and \(\sigma(x) \equiv \Coh \Theta B \rho\). Then \((\tau \circ \sigma) (x) \equiv \Coh \Theta B {\tau \circ \rho}\) and so \(\sub t \tau\) admits an insertion. We must have \(t' \equiv \Coh {\insertion \Delta x \Theta} {\sub B {\kappa}} {\insertion \sigma x \rho}\) and \(\sub t \tau\) reduces by insertion to \(\Coh {\insertion \Delta x \Theta} {\sub B {\kappa}} {\insertion {(\tau \circ \sigma)} x {(\tau \circ \rho)}}\). It remains to show that \(\tau \circ {(\insertion \sigma x \rho)} \equiv \insertion {(\tau \circ \sigma)} x {(\tau \circ \rho)}\) but this is given by \cref{lem:inserted-sub}.
    \end{itemize}
  \end{itemize}
\end{proof}

We end this section by showing that reduction agrees with equality in the following way.

\begin{definition}
  Let \(\rto'\) be the reduction \(\rto\) restricted to terms that are well typed, and \(\rto'_{\mathsf{rts}}\) be the reflexive transitive symmetric closure of this reduction.
\end{definition}

\begin{prop}\label{prop:reduct-equality}
  On valid terms, types and substitutions, the equality relation \(=\) agrees with \(\rto'_{\mathsf {rts}}\).
\end{prop}
\begin{proof}
  It is clear that if \(p \rto'_{\mathsf{rts}} q\) then \(p = q\). Therefore we only need to prove the other direction. As in \cite{finster20_type_theor_stric_unital_categ}, we define a term \(p\) to be \emph{conservative} if for all \(q\) with \(p = q\), we have \(p \rto'_{\mathsf{rts}} q\). The only case that doesn't follow almost immediately from inductive hypothesis is the case where \(p = q\) is by insertion. Therefore suppose \(p \equiv \Coh \Delta A \sigma\) and \(\sigma(x) = \Coh \Theta U \tau\) and \(q \equiv \Coh {\insertion \Delta x \Theta} {\sub A \iota} {\insertion \sigma x \tau}\). Then by induction on subterms we have that \(\sigma(x) \rto'_{\mathsf{rts}} \Coh \Theta U \tau\) and so \(p \rto'_{\mathsf{rts}} \Coh \Delta A {\sigma'} \rto' q\) where \(\sigma'\) is \(\sigma\) with \(\sigma(x)\) replaced with \(\Coh \Theta U \tau\).
\end{proof}

\subsection{Termination}
\label{sec:termination}

In this section we prove termination of the reduction scheme presented above. To do this we will map each term to an ordinal, show that all reductions reduce this ordinal.

\begin{definition}
  The natural sum of ordinals \(\alpha \+ \beta\) is defined by simultaneous induction on \(\alpha\) and \(\beta\) by:
  \begin{alignat*}{2}
    &0 \+ 0 & &= 0\\
    &\alpha \+ \beta &&= \sup _{
      \substack{\alpha' < \alpha\\ \beta' < \beta}
    }
    \{ S(\alpha \+ \beta'), S(\alpha' \+ \beta) \}
  \end{alignat*}
\end{definition}

\begin{lemma}\label{lem:natsum}
  The natural sum of ordinals is associative, commutative, and strictly monotone in both arguments.
\end{lemma}
\begin{proof}
  See~\cite{lipparini16_infin_natur_sum}.
\end{proof}

\begin{definition}
  Define the syntactic depth \(\sd\) of substitutions, typed, and terms by mutual induction:
  \begin{itemize}
  \item For substitutions we have:
    \[\sd(\langle t_0, \dots, t_n \rangle) = \bighash_{i=0}^n t_i\]
  \item For types we have \(\sd(\star) = 0\) and:
    \[\sd(\arr s A t) = \sd(s) \+ \sd(A) \+ \sd(t)\]
  \item For variables \(x\), \(\sd(x) = 0\).
  \item For coherences define:
    \[\sd(\Coh \Delta A \sigma) = \omega^{\dim(A)} \+ \sd(A) \+ \sd(\sigma)\]
  \end{itemize}
\end{definition}


\begin{lemma}\label{lem:sd-decrease}
  For terms \(t \rto t'\), \(\sd(t) > \sd(t')\), for types \(A \rto A'\), \(\sd(A) > \sd(A')\), and for substitutions \(\sigma \rto \sigma'\), \(\sd(\sigma) > \sd(\sigma')\).
\end{lemma}
\begin{proof}
  By mutual induction on term reductions \(t \rto t'\), type reductions \(A \rto A'\) and substitution reductions \(\sigma \rto \sigma'\):
  \begin{itemize}
  \item For substitutions, if \(\sigma \rto \sigma'\) then suppose \(\sigma = \langle t_0, \dots, t_n \rangle\) and \(\sigma' = \langle t_0', \dots, t_n' \rangle\). Then for some \(i\) we must have \(t_i \rto t_i'\) and \(t_k = t_k'\) for \(k \neq i\). By inductive hypothesis \(\sd(t_i) > \sd(t_i')\) and so \(\sd(\sigma) > \sd(\sigma')\).
  \item For types, \(\star\) admits no reductions so suppose \(A \equiv \arr s B t\). If \(A' \equiv \arr {s'} B t\) where \(s \rto s'\) then by inductive hypothesis \(\sd(s) > \sd(s')\) and so \(\sd(A) > \sd(A')\). The cases where \(A' \equiv \arr s {B'} t\) and \(A' \equiv \arr s B {t'}\) follow similarly.
  \item For term reductions \(t \rto t'\) we first split on the reduction. Note that as variables do not reduce, we must have that \(t \equiv \Coh \Delta A \sigma\).
    \begin{itemize}
  \item Argument reduction: If \(t \rto t'\) is an argument reduction then \(t' \equiv \Coh \Delta A {\sigma'}\) where \(\sigma \rto \sigma'\). By inductive hypothesis, \(\sd(\sigma) > \sd(\sigma')\) and so \(sd(t) > \sd(t')\).
  \item Cell reduction: If \(t \rto t'\) is a cell reduction then we have \(t' \equiv \Coh \Delta {A'} \sigma\) with \(A \rto A'\). By inductive hypothesis \(\sd(A) > \sd(A')\) and so \(\sd(t) \geq \sd(t')\).
  \item Insertion: If \(t \rto t'\) is an insertion along a variable \(x\) then let \(\sigma(x) \equiv \Coh \Theta B \tau\). Now let \(\mu \equiv \insertion \sigma x \tau\) and let \(I \equiv \insertion \Delta x \Theta\). Then:
    \begin{align*}
      \omega^{\dim(A)} \+ \sd(\mu) &= \omega^{\dim(A)} \+ \bighash_{y \in \Var(I)} \sd(\mu(y))\\
               &= \omega^{\dim(A)} \+ \bighash_{y \in \Var(I) \cap \Var(\Theta)} \sd(\mu(y)) \+ \bighash_{y \in \Var(I) \cap \Var(\Delta)} \sd(\mu(y)) \\
               &= \omega^{\dim(A)} \+ \bighash_{y \in \Var(\Theta)} \sd(\tau(y)) \+ \bighash_{y \in \Var(I) \cap \Var(\Delta)} \sd(\sigma(y))\\
               &= \omega^{\dim(A)} \+ \sd(\tau) \+ \bighash_{y \in \Var(I) \cap \Var(\Delta)} \sd(\sigma(y))\\
               &\geq \sd(\sigma(x)) \+ \bighash_{y \in \Var(I) \cap \Var(\Delta)} \sd(\sigma(y))\\
               &= \bighash_{y \in (\Var(I) \cup \{x\}) \cap \Var(\Delta)} \sd(\sigma(x))\\
               &\leq \bighash_{y \in \Var(\Delta)} \sd(\sigma(y))\\
               &= \sd(\sigma)
    \end{align*}
    Then by assumption, \(t\) is not a coherence on a disc, and so is not an identity. Therefore:
    \begin{align*}
      \sd(t) &=  \omega^{\dim(t)} \+ \sd(A) \+ \sd(\sigma)\\
             &\geq \omega^{\dim(t)} \+ \sd(\sigma) \\
             &\geq \omega^{\dim(t)} \+ \omega^{\dim(t)} \+ \sd(\insertion \sigma x \tau)\\
             &> \omega^{\dim(t)} \+ \sd(\sub A \kappa) \+ \sd((\insertion \sigma x \tau)) \\
             &= \sd(t')
    \end{align*}
    where \(\omega^{\dim(t)} > \sd(\sub A \kappa)\) due to dimension.
  \end{itemize}
  \end{itemize}

  This covers all cases, and so completes the proof.
\end{proof}

\begin{theorem}\label{thm:red-terminates}
  The reduction \(\rto\) is strongly normalising on types, terms, and substitutions.
\end{theorem}
\begin{proof}
  Suppose an infinite reduction sequence \(t_0 \rto t_1 \rto \dots\) exists. Then by \cref{lem:sd-decrease}, \(sd(t_0) > \sd(t_1) > \dots\) is an infinitely decreasing sequence of ordinals, which cannot exist.
\end{proof}

\subsection{Confluence}
\label{sec:confluence}

In this section we show that \(\rto'\) is locally confluent, which implies along with \cref{thm:red-terminates} that \(\rto'\) is globally confluent and hence complete.

\begin{definition}
  Mutually define a term \(t\) being \emph{regular} and its \emph{regular height}, an extended natural number as follows:
  \begin{itemize}
  \item If \(t\) is a variable then it is regular and its regular height is \(\infty\)
  \item If \(t \equiv \Coh \Delta U \sigma\), then \(t\) is regular if and only if \(\Delta\) is not a disk, \(U\) is the unbiased type over \(\Delta\), every term in \(\sigma\) is regular, and for all \(x \in \LM(\Delta)\), the branching height of \(x\) is less than the regular height of \(\sigma(x)\). If \(t\) is regular, then its regular height is given by the linear height of \(\Delta\).
  \end{itemize}
\end{definition}

If \(I,\iota_i,\iota_e\) is an insertion of \((\Theta,B)\) into \((\Delta,x)\), with \(U\) is the unbiased type over \(\Delta\) and \(B\) is the unbiased type over \(\Theta\), then all terms contained in \(\sub U {\iota_e}\) are regular, as all the terms in \(U\) are unbiased composites of variables, and we substitute in unbiased composites which satisfy the height condition from the conditions of insertion.

Further note that if the support of a regular term is a disc, then it must be a variable, as any unbiased composite will have multiple locally maximal variables in its support.

\begin{lemma}
  If \(t\) is a regular term over a pasting diagram, then it reduces to the unbiased composite of its support.
\end{lemma}
\begin{proof}
  Proceed by mutual reduction on dimension and subterms. The hypothesis on dimension tells us that if a regular term undergoes an insertion, then it reduces to another regular term, (in other words the type reduces to the unbiased type). Further, we must have that the regular height of the new reduced term is at least the regular height of the original term. Therefore take an arbitrary term \(t\), and repeatedly perform insertions followed by making everything unbiased. By termination, we will eventually get term \(t'\) with \(t \rto t'\) and \(t'\) having no more possible insertions.

  If \(t'\) is a variable then we are done. Otherwise, \(t'\) is an unbiased composition. By construction, \(t'\) does not admit any insertions, so all locally maximal dimensional arguments must be variables (or \(t'\) would not be regular). Now suppose there is an argument \(x\) is not a variable. Then there is such an argument \(u\) with maximal dimension. This can not be locally maximal so it is the source or target of a variable. Then its support must be a disc, but this contradicts its regularity.
\end{proof}

\begin{cor}\label{cor:red-to-unbiased}
  If \(t\) is an unbiased composite and admits an insertion, then the resulting term reduces to an unbiased composite.
\end{cor}

\begin{theorem}
  The reduction \(\rto'\) is locally confluent. In other words if \(a\), \(b\), and \(c\) are valid terms, with \(a \rto b\) and \(a \rto c\) then there exists valid \(u\) with \(b \rto_{\mathsf{rt}} u\) and \(c \rto_{\mathsf{rt}} u\).
\end{theorem}
\begin{proof}
  Proceed by mutual induction on subterms and dimension. By induction on dimension we have that \(\rto\) is locally confluent on \(m\)-dimensional terms for \(m < \dim(a)\) and so is globally confluent (for dimension \(m\)). Further, as \(\rto\) is terminating (by \cref{thm:red-terminates}) and agrees with equality (by \cref{prop:reduct-equality}), we have that any equal terms with dimension less than \(\dim(a)\) have a common normal form. It also holds that if \(s =_m t\) for some \(m < \dim(a)\), then there is \(u\) with \(s \rto_{\mathsf{rt}} u\) and \(t \rto_{\mathsf{rt}} u\).

  Suppose \(a \rto b\) and \(a \rto c\). Split into cases.
  \begin{itemize}
  \item Both are insertions. If they are insertions on the same variable, then \(b \equiv c\) and we are done. Otherwise let \(a \equiv \Coh \Gamma A \sigma\), and suppose \(b\) is an insertion along \(x\) with \(\sigma(x) \equiv \Coh \Delta B \tau\) and \(c\) is an insertion along \(y\) with \(\sigma(y) \equiv \Coh \Theta C \mu\). Further we must have that \(x\) and \(y\) are locally maximal, so in particular are not in each other's support.

\[\begin{tikzcd}
        & {D^n} & \Delta \\
        {D^m} & \Gamma & {\insertion \Gamma x \Delta} \\
        \Theta
        \arrow["\overline{x}", from=1-2, to=2-2]
        \arrow["{\overline{\Coh \Delta B \id}}", from=1-2, to=1-3]
        \arrow[from=2-2, to=2-3]
        \arrow[from=1-3, to=2-3]
        \arrow["\overline{y}", from=2-1, to=2-2]
        \arrow["{\overline{\Coh \Theta C \id}}"', from=2-1, to=3-1]
        \arrow["\lrcorner"{anchor=center, pos=0.125, rotate=180}, draw=none, from=2-3, to=1-2]
\end{tikzcd}\]
In the diagram above we construct the insertion along \(x\). Now we can notice that the composite \(D^m \overset{y}{\to} \Gamma \to \insertion \Gamma x \Delta\) is just the substitution \(\overline{y}\), and so we can perform another insertion (as the branching height of \(y\) must be the same in \(\Gamma\) and \(\insertion \Gamma x \Delta\)). We then get the following diagram:
\[\begin{tikzcd}
        & {D^n} & \Delta \\
        {D^m} & \Gamma & {\insertion \Gamma x \Delta} \\
        \\
        \Theta && {\insertion {(\insertion \Gamma x \Delta)} y \Theta}
        \arrow["\overline{x}", from=1-2, to=2-2]
        \arrow["{\overline{\Coh \Delta B \id}}", from=1-2, to=1-3]
        \arrow[from=2-2, to=2-3]
        \arrow[from=1-3, to=2-3]
        \arrow["\overline{y}", from=2-1, to=2-2]
        \arrow["{\overline{\Coh \Theta C \id}}"', from=2-1, to=4-1]
        \arrow["\lrcorner"{anchor=center, pos=0.125, rotate=180}, draw=none, from=2-3, to=1-2]
        \arrow[from=2-3, to=4-3]
        \arrow[from=4-1, to=4-3]
        \arrow["\lrcorner"{anchor=center, pos=0.125, rotate=180}, draw=none, from=4-3, to=2-1]
      \end{tikzcd}\]
    and therefore \(\insertion {(\insertion \Gamma x \Delta)} y \Theta\) is the colimit of:
\[\begin{tikzcd}
        & {D^n} & \Delta \\
        {D^m} & \Gamma \\
        \Theta
        \arrow["\overline{x}", from=1-2, to=2-2]
        \arrow["{\overline{\Coh \Delta B \id}}", from=1-2, to=1-3]
        \arrow["\overline{y}", from=2-1, to=2-2]
        \arrow["{\overline{\Coh \Theta C \id}}"', from=2-1, to=3-1]
      \end{tikzcd}\]
    Similarly, \(\insertion {(\insertion \Gamma y \Theta)} x \Delta\) is a colimit of the same diagram and so is equal. Further, it can easily be checked that \(\insertion {(\insertion \sigma y \mu)} x \tau\) and \(\insertion {(\insertion \sigma x \tau)} y \mu\) both satisfy the necessary equations to be equal by uniqueness of the map out of a colimit. It is then simple to check that the maximal dimension parts of these two substitutions are actually definitionally equal by case analysis and so we get
    \begin{align*}
      b &\rto \Coh {\insertion {(\insertion \Gamma x \Delta)} y \Theta} {\sub A {\kappa_{(\insertion \Gamma x \Delta), y , \Theta} \circ \kappa_{\Delta, x, \Gamma}}} {\insertion {(\insertion \sigma x \tau)} y \mu} \\
        &=_{\dim(a)} \Coh {\insertion {(\insertion \Gamma y \Theta)} x \Delta} {\sub A {\kappa_{(\insertion \Gamma y \Theta), x, \Delta} \circ \kappa_{\Theta, y, \Gamma}}} {\insertion {(\insertion \sigma y \mu)} x \tau}\\
      &\leftsquigarrow c
  \end{align*}
    which gives us a way to construct a suitable \(u\) by the comment at the start of the proof.
  \item Now suppose both reductions are argument reductions. If these occur on different arguments, we simply apply both to get a suitable \(u\), and if they are the same argument then we use the inductive hypothesis on subterms.
  \item Now suppose \(a \rto b\) is an insertion and \(a \rto c\) is an argument reduction. The only case which is hard is when the insertion and argument reduction are along the same variable. Therefore we can assume we are in the case where \(a \equiv \Coh \Gamma A \sigma\) and \(\sigma(x) \equiv \Coh \Delta B \tau\) and \(\tau(y) \equiv \Coh \Theta C \mu\).

    Firstly assume that \(\Delta\) is a disc context. Then \(\sigma(x)\) just reduces to \(\tau(y)\) by performing the argument insertion first. If instead the head insertion is performed first then the head pasting diagram is unchanged (up to variable renaming) but the arguments can be changed. We see that \(x\) must be sent to \(\tau(y)\). Therefore \(b =_{\dim(x)} c\) and so we can find suitable \(u\) as before.

    Now we assume \(\Delta\) is not a disc context. First consider the case where the insertion is done first. Then we get the following diagram:

\[\begin{tikzcd}
        & {D^n} & \Gamma \\
        {D^m} & \Delta & {\insertion \Gamma x \Delta} \\
        \\
        \Theta && {\insertion {(\insertion \Gamma x \Delta)} y \Theta}
        \arrow["\overline{x}", from=1-2, to=1-3]
        \arrow["\overline{y}", from=2-1, to=2-2]
        \arrow[from=2-2, to=2-3]
        \arrow[from=1-3, to=2-3]
        \arrow["{\overline{\Coh \Delta B \id}}"', from=1-2, to=2-2]
        \arrow["{\overline{\Coh \Theta C \id}}"', from=2-1, to=4-1]
        \arrow["\lrcorner"{anchor=center, pos=0.125, rotate=180}, draw=none, from=2-3, to=1-2]
        \arrow[from=4-1, to=4-3]
        \arrow[from=2-3, to=4-3]
        \arrow["\lrcorner"{anchor=center, pos=0.125, rotate=180}, draw=none, from=4-3, to=2-1]
      \end{tikzcd}\]
    where we can perform the second insertion as the composite \(D^m \overset y \to \Delta \to \insertion \Gamma x \Delta\) is just the variable \(y\) (as \(\Delta \to \insertion \Gamma x \Delta\) is variable to variable) and \(y\) has the same branching height in \(\Delta\) and \(\insertion \Gamma x \Delta\) as \(\Delta\) is not a disc context.

    Now consider first applying the argument insertion to get this diagram:
\[\begin{tikzcd}
        & {D^n} & \Gamma \\
        {D^m} & \Delta \\
        \Theta & {\insertion \Delta y \Theta}
        \arrow["\overline{x}", from=1-2, to=1-3]
        \arrow["\overline{y}", from=2-1, to=2-2]
        \arrow["{\overline{\Coh \Delta B \id}}"', from=1-2, to=2-2]
        \arrow["{\overline{\Coh \Theta C \id}}"', from=2-1, to=3-1]
        \arrow[from=2-2, to=3-2]
        \arrow[from=3-1, to=3-2]
        \arrow["\lrcorner"{anchor=center, pos=0.125, rotate=180}, draw=none, from=3-2, to=2-1]
      \end{tikzcd}\]

    Now note that there is a second map from \(D^n\) to \(\insertion \Delta y \Theta\) given by \(\Coh {\insertion \Delta y \Theta} U \id\) where \(U\) is the unbiased type of \(\insertion \Delta y \Theta\) and by an insertion and \cref{cor:red-to-unbiased} these are equal. As \(\Delta\) is not a disc context, the linear height of \(\insertion \Delta y \Theta\) is at least as high as \(\Delta\) and so we can complete the second insertion to get:
\[\begin{tikzcd}
        & {D^n} && \Gamma \\
        {D^m} & \Delta \\
        \Theta & {\insertion \Delta y \Theta} && {\insertion \Gamma x {(\insertion \Delta y \Theta)}}
        \arrow["\overline{x}", from=1-2, to=1-4]
        \arrow["\overline{y}", from=2-1, to=2-2]
        \arrow["{\overline{\Coh \Delta B \id}}"', from=1-2, to=2-2]
        \arrow["{\overline{\Coh \Theta C \id}}"', from=2-1, to=3-1]
        \arrow[from=2-2, to=3-2]
        \arrow[from=3-1, to=3-2]
        \arrow["\lrcorner"{anchor=center, pos=0.125, rotate=180}, draw=none, from=3-2, to=2-1]
        \arrow[""{name=0, anchor=center, inner sep=0}, "{\overline{\Coh {\insertion \Delta y \Theta} U \id}}", curve={height=-30pt}, from=1-2, to=3-2]
        \arrow[from=1-4, to=3-4]
        \arrow[from=3-2, to=3-4]
        \arrow["{=}"{description}, Rightarrow, draw=none, from=2-2, to=0]
\end{tikzcd}\]
Then as before, both of these insertions arise from the same colimit and we can apply the same reasoning as the insertion-insertion case, which completes the proof.
  \end{itemize}
\end{proof}

\section{Typechecking}
\label{sec:typechecking}

Using the reduction in the previous section, we can show that equality and typechecking are decidable.

\begin{definition}
  Define the normalisation function \(N\) as a function taking terms to terms, types to types, and substitutions to substitutions, given by repeatedly applying the innermost leftmost reduction until a normal form is found.
\end{definition}

\begin{prop}
  If \(s\) and \(t\) are valid terms, then \(s = t\) if and only if \(N(s) \equiv N(t)\). A similar result holds for valid typed \(A\) and \(B\), and valid substitutions \(\sigma\) and \(\tau\).
\end{prop}

This gives us a way to decide equality between terms which are valid, as we can simply compute the normalisation of each term and checking syntactic equality is clearly decidable. We can now show how this can be used to perform decidable typechecking.

\begin{theorem}
  Typechecking on terms, types, substitutions, and contexts is decidable.
\end{theorem}
\begin{proof}
  We work by mutual induction on terms, types, substitutions, and contexts.
  \begin{itemize}
  \item For typechecking whether \(\Gamma \vdash t : A\), we split into cases:
    \begin{itemize}
    \item If \(t\) is a variable, then we check whether \(t \in \Gamma\). If it is not, then \(t\) is not valid in \(\Gamma\). If \(t : B\) in \(\Gamma\) then we just need to check if \(A = B\).
    \item If \(t \equiv \Coh \Delta B \sigma\), then for \(t\) to be valid, it is necessary that \(\Gamma \vdash \sigma : \Delta\), \(\Delta \vdash_p\), and \(\Delta \vdash B\). It can be easily checked that \(\Delta\) is a pasting diagram, and the other two can be checked by inductive hypothesis. We then need to further check that the support conditions hold for either the coherence or composition rule. Finally, if all other conditions hold then \(\Gamma \vdash t : A\) if and only if \(A = \sub B \sigma\).
    \end{itemize}
  \item For typechecking \(\Gamma \vdash A\), if \(A \equiv \star\) then \(A\) is valid, and if \(A \equiv \arr s B t\), then \(A\) is valid if and only if \(\Gamma \vdash B\), \(\Gamma \vdash s : B\), and \(\Gamma \vdash t : B\), all of which can be checked by inductive hypothesis.
  \item Typechecking for substitutions similarly follows from inductive hypothesis.
  \item The empty context is always valid, and \(\Gamma, x : A\) can be checked by checking whether \(\Gamma \vdash\) and \(\Gamma \vdash A\).
  \end{itemize}
  This gives us a fully decidable procedure for typechecking. It can also be noted that we can extract a type inference algorithm for terms from this proof.
\end{proof}

\todo{Talk about implementation?}

\printbibliography{}
\end{document}